\documentclass[12pt]{article}
\usepackage{epsfig,psfrag,amsmath,amssymb,latexsym}
\usepackage{amscd }
\usepackage{color}
\usepackage{amsfonts}
\usepackage{graphicx}

\usepackage[pdftex,colorlinks,linkcolor=blue,menucolor=red,bookmarks=true]{hyperref}

\usepackage{makeidx,amsfonts, amsthm, amssymb, amsmath, graphicx, bbm, color, xcolor,  wrapfig, caption,subcaption, listings, float, mathtools,  url, dsfont}
\usepackage{epstopdf}

\usepackage[top=1in, bottom=1in,left=1in,right=1in]{geometry}

\newtheorem{theorem}{Theorem}[section]

\newtheorem{example}{Example}[section]

\numberwithin{equation}{section}
\begin{document}
\title{Recovery of spectrum from estimated covariance matrices
and statistical kernels for machine learning and big data} 
\author{Saba Amsalu, \thanks{University of Addis Abeba, Ethiopia}
	\and	
		Juntao Duan 
		\thanks{School of Mathematics, Georgia
		Institute of Technology, Atlanta, GA 30332, 
		\hfill\break Email: jt.duan@gatech.edu}
	\and Heinrich Matzinger 
		\thanks{School of Mathematics, Georgia
		Institute of Technology, Atlanta, GA 30332
		\hfill\break Email: matzi@math.gatech.edu}
	\and Ionel Popescu 
		\thanks{School of Mathematics, Georgia
		Institute of Technology, Atlanta, GA 30332, and 
		IMAR and FMI, Bucharest Romania.
		\hfill\break Email: ipopescu@math.gatech.edu} 
		\thanks{
		I.P. acknowledges that this work was partially supported by 
	Simon Collaboration grant no. 318929 and  
	UEFISCDI grant PN-III-P4-ID-PCE-2016-0372.} 	
      }

\maketitle
\begin{abstract}
Let $X$ be a $n\times p$ random matrix with independent rows
all distributed as the random vector $\vec{X}=(X_1,X_2,\ldots,X_p)$
and having $0$ expectation. The sample covariance matrix
is then defined by $\widehat{COV}[\vec{X}]=X^tX/n$ and is
the maximum likelihood  estimate of the covariance matrix
 $COV[\vec{X}]$ when the data is normal.  We are interested in the regime where
 $n$ and $p$ both go to infinity at the same time
so that $n/p$ converges to a non-zero number $\alpha>0$. In that case,
the spectrum of the sample covariance converges to the free product of $\mu$
with a Marchenko-Pastur distribution \cite{pastur1967}, provided the spectrum
of the covariance converges to $\mu$. In theory, from the limiting 
distribution of the estimated covariance matrix, one could retrieve $\mu$ using the 
$S$-transform from free probability. 

In many practical applications, the spectrum
of the ``actual'' covariance contains essential information about
the structure of the data at hand.  There is a family
 of spectrum-reconstruction-algorithms,
 which attempt to discretize and adapt the free probability
infinite-dimensional-recovery. This approach, pioneered by  EI Karoui 
\cite{elkaroui2008}, was further developed by Silverstein \cite{silverstein1995}, Bai and Yin \cite{bai1988}, Yin, Bai and Krishnaiah \cite{yin1983}
and recently Ledoit and Wolf  \cite{ledoit2012}, \cite{li2013} and \cite{ledoit2015}.
There is also the shrinkage method pioneered by Stein \cite{Stein}, 
see also Bickel and Levina \cite{Bickel} and Donoho \cite{Donoho}. 
Another type of approach is based on the moments of the spectral distributions see for example Kong and Valliant \cite{Valiant}. This can be very useful as parallel method. Finally, there are also
Physicists Burda, G\"orlich and Jarosz, working on this problem \cite{Burda}.

We present a simple algebraic formula  \eqref{library}
for reconstructing the spectrum of the covariance matrix
given only the sample covariance spectrum. For normal data
with $n/p\geq 2$ and the spectrum taken from real data, this approximation
allows to recover the true underlying spectrum quite closely (See Figure 
\ref{eigenvectormethod}). For spectrum which is too flat however 
it does not work well (See Figure \ref{constantspectrum}).

We also introduce a second algorithm which is a  fixed point method, 
with remarkable reconstruction of the spectrum when $n/p\ge 1$ in both simulated and 
real data. For this, we take real life data with $n$ much bigger than $p$ to estimated 
the ``true covariance'''s spectrum. Then, we compare it to our reconstruction
using a small subset of the real-life data for which $n/p$ is much smaller. 
(See Figure \ref{realdata}, where $n/p=1$).

Finally we present a proof of concept for big data, by showing that
our method applied to a reduced real life data set, allows to recover
almost perfectly the underlying spectrum of the full-data set. (See 
Figure \ref{Reconstructed_by_scaling_down}). 
Our method does not seem to work well on sparse data, 
though this is not our main focus in this note. 

Our first approximation  is derived  in a purely
probabilistic way using the formula for
the spectral norm of error matrix \cite{Kolt1}
 by Koltchinskii and Lounici.  
  
\end{abstract}
\tableofcontents

\section{Introduction}
\label{intro}

\subsection{The setup}

We consider a random matrix $Z$ of dimension $n\times p$ which represents some data.
We will be primarily interested in the eigenvectors of the expected matrix
$E[Z^tZ]$ called {\it Principal Components}. The principal components
corresponding to the highest eigenvalues, often carry some structural
information about the underlying data. Instead
of being given the expected matrix, only $Z^tZ$ is usually observable.
We will take the eigenvectors of $Z^tZ$ and hope they approximate
the eigenvectors of $E[Z^tZ]$ well. This is the setting for covariance matrix
analysis via principal components, as well as analysis of 
statistical kernels via principal components.

For example, when the rows of $Z$ have $0$ expectation
and  are i.i.d., the expectation of $Z^tZ/n$ corresponds to the covariance matrix
\begin{equation}
\label{covariance}
COV(\vec{Z})=E[\frac{Z^tZ}{n}].
\end{equation}
where $\vec{Z}$ has the distribution of a row of $Z$.
In this case of i.i.d. rows with $0$ expectation, 
the sample covariance matrix is  given
\begin{equation}
\label{estimated}
\widehat{COV}[\vec{Z}]:=\frac{Z^tZ}{n}.
\end{equation}
It is an unbiased maximum-likelihood estimator of the 
true covariance matrix when the data is normal.

In this paper, we present two new algorithms
which allow in many practical situations
to get very close to retrieving the spectrum of $E[Z^tZ]$
given only the spectrum of $Z^tZ$. Our algorithms  to our knowledge
allow a level of precision with most real-life data which is remarkable. 
Our reconstruction of the real spectrum is often  precise enough to determine
details of the spectrum of $E[Z^tZ]$ like for example ``elbows''.

Let $X$ be a $n \times p$ matrix 
with i.i.d. multivariate standard normal random variables. 
The sample covariance matrix is then defined by 
$\widehat{COV}[{X}]=\frac{1}{n} X^tX$ and is
a maximum-likelihood estimate of the covariance matrix
$COV[{X}]$, which is identity. We are interested in the 
high dimensional setting,  where $n$ and $p$ both go 
to infinity at the same time so that $\frac{p}{n}$ 
converges to a non-zero fixed limit. In 1967, Vladimir 
Mar{\v c}enko and Leonid Pastur \cite{pastur1967} 
successfully constructed the limiting law, which is now named after the authors. 

For the case of  $COV[{X}]$ is not identity, we assume 
its spectrum admits a limiting law $F^{\Sigma}$. Then sample 
spectrum is a free product of Mar{\v c}enko-Pastur law and $F^{\Sigma}$. 
By computing the $S$-transform explicitly, one can obtain a 
formula of the limiting law of the sample spectrum as for instance discussed in 
Bai and Yin \cite{bai1988}, Yin, Bai and Krishnaiah \cite{yin1983}, 
Silverstein \cite{silverstein1995}, and many others.  The main result is 
summarized as follows.

\begin{theorem}
	Assume the following. 
	\begin{enumerate}
		\item{} The entries of $X_p=(X_{i,j})_{n\times p}$ are 
			i.i.d. real random variables for all $p$. 
		
		\item{} $E[X_{1,1}]=0$, $E[|X_{1,1}|^2]=1$. 
		
		\item{} Let $p/n \to c >0$ as $p\to \infty$. 
		
		\item{} Let $\Sigma_p$ $(p\times p)$ be non-negative definite symmetric random matrix with spectrum distribution $F^{\Sigma_p}$ (If $\{\lambda_i\}_{1\le i \le p}$ are the eigenvalues of $\Sigma_p$, then $F^{\Sigma_p}=\sum_{1}^{p} \frac{1}{p} \delta_{\lambda_i}(x)$)  such that $F^{\Sigma_p}$ almost surely converges weakly to $F^{\Sigma}$ on $[0, \infty)$. 
		
		\item{} $X_p$ and ${\Sigma}_p$ are independent.
	\end{enumerate}
	Then the spectrum distribution of $W_p= \frac{1}{n}{\Sigma}_p^{1/2}X_p^T X_p {\Sigma}_p^{1/2}$, denoted as $F^{W_p}$  almost surely converges weakly to $F^W$. $F^W$ is the unique probability measure whose Stieltjes transform $m(z)= \int \frac{d F^W(x)}{x-z}$, $z\in \mathbb{C}^+$ satisfies the equation
	\begin{equation} \label{eqn:sample spectrum}
	-\frac{1}{m}=z- c\int \frac{t }{1+tm} d F^{\Sigma}(t) \quad \forall z \in \mathbb{C}^+
	\end{equation}
	
\end{theorem}

In many real life application,  the  spectrum of the true covariance
matrix contains essential information about the underlying structure
of the data. Now, this is not the spectrum of the estimated
covariance if there is a big difference.  
For example if $n=p/2$, then the estimated covariance matrix 
is defective and has at least half of its eigenvalues equal to $0$.  
In financial analysis, a concrete example is the daily stock returns. 
The spectrum of the underlying covariance matrix are bounded away from $0$ 
by a term of order $O(1)$, since stocks are not linearly dependent. 
Thus spectrum of the defective sample covariance matrix ``makes'' 
an error of order at least $O(1)$. 

One important approach is using the free probability results 
and attempts at solving the equation \eqref{eqn:sample spectrum} 
to get an estimator of the true spectrum $F^{\Sigma}$. 
First such work is due to EI Karoui \cite{elkaroui2008}, 
and then Bai etc. \cite{bai2010estimation}, and recently by 
Ledoit and Wolf  \cite{ledoit2012}, \cite{li2013} and \cite{ledoit2015}. 
It's not surprising that as dimensions grow, consistency is achieved 
in theory, by the free probability approach. We say ``in theory'',
because for real-life-data there are usually some eigenvalues of the spectrum
which are of order $O(p)$ and others of order $O(1)$,
so that the spectrum of the underlying covariance matrix does
not converge to a limit. Furthermore, a disadvantage 
of the free-probability approach is that 
the recovered spectrum can still be far from the true spectrum for 
small or moderate size of $p$.

\subsection{Our results}
 
In the current paper, we first  present a simple algebraic formula  
\eqref{approxlambda} for the sample eigenvalues which approximates 
the true covariance eigenvalues. For normal data, this approximation 
is quite efficient and gives quite good corrections
as soon as $n/p\ge2$ and the true spectrum is not too
flat (see for instance Figure \ref{eigenvectormethod}). 
In addition to this we present a more precise recovery, 
we introduce a fixed point iteration method. 
This method uses the matrix of the eigenvectors of  the true covariance
matrix expressed in the basis of the eigenvectors of the sample covariance.

The second algorithm  is more costly but can reconstruct 
almost perfectly if the true spectrum is not too wild. 
Our methods are also applied to real data which turns out to 
work very well as long as the data at hand is not sparse. 

Our formula \eqref{library} leading to first approximation, looks similar 
to the $S$-transform from free probability and thus resembles the free probability
approach. Nonetheless, we derived it in a purely probabilistic way using
the formula for the size of the  norm of the error matrix of covariance
\cite{Kolt1} 
proven by Koltchinskii and Lounici. Our method is the only one which, 
to our knowledge, is often precise enough to allow to detect important 
features such as ``elbows'' in the spectrum of real data.
Our formula \eqref{library} for approximating the difference 
in spectrum between sample covariance and true covariance, could be written
with the distribution function in a similar way to
\eqref{eqn:sample spectrum} as follows
\begin{equation}\label{our}
-\int\frac{y}{y-invF^{W_p}(x)}dF^{W_p}(y)\approx
\frac{invF^{W_p}(x)}{invF^{\Sigma_p}(x)}-1,
\end{equation}
for $x\in[0,1]$. Here $invF(.)$ designates the inverse 
of the distribution function $F(.)$.  Also, when $p$ is finite,
then the integral on the right side of \eqref{our} is an integration with 
respect to a discrete measure, and one needs to leave out
$y=inv F^{W_p}(x)$ (and even $3$ to $5$ $y$'s closest to 
$inv F^{W_p}(x)$) for the formula to work. At the limit,
the right side of \eqref{our} is an improper integral.
Now, the difference between the free probability approach and ours
is that though for many real life situations, our
formula \eqref{our} is very precise, we don't expect it to converge
to an exact solution as $p$ goes to infinity.
In other words, as $p$ goes to infinity, \eqref{our}
should not become an exact equation, but should remain an approximation.

In Subsection \eqref{aymptoticmethod}, we show how 
the approximation  \eqref{library} is derived, without
providing a rigorous bound for the error of that approximation. 
We just note that for practical purposes, the approximation \eqref{library}
seems quite good.

Our second spectrum reconstruction algorithm is a fixed point
method. We present it in Subsection \ref{SectionEigen}.
There we introduce  a  map \eqref{themap}  defined with the help
of the sample covariance spectrum and  related to the matrix of the second moments of the coefficients
for the eigenvectors. That is, the eigenvectors of the true covariance matrix
expressed in the basis of the eigenvectors of the sample covariance matrix.
The map \eqref{themap} can be applied to any spectrum, but is defined
with the help of the sample covariance. In Subsection \ref{SectionEigen}, 
we show that for this map and up to a smaller order term, the
spectrum of the true covariance matrix is a fixed point. We prove nothing
else, but for practitioners this is already very valuable. Indeed,
a continuous map from a compact space into itself could be chaotic,
periodic or have fixed points. We know from practice that our map
is not chaotic nor periodic, so it must have fixed points. Hopefully, 
in the terminology of dynamical systems, this fixed point is an attractor 
and this would explain why we have this great convergence, however we 
do not have a proof of this results yet.  
Hence, the main problem would be that we converge to a wrong fixed point,
(i.e. not to the desired spectrum of the true covariance). 
In practice this is not happening, except if the true spectrum is very flat
or if $n/p$ is quite smaller than say about $0.5$. On top of it,
we can use our first method to find a good enough starting point
for the eigenvalues fixed point method. The fixed point algorithm presented
in Subsection \ref{SectionEigen} consists then simply in applying
the map \eqref{themap} until it stabilizes. In many examples, one iteration 
of \eqref{themap} is enough to get a good approximation of the spectrum
of the true covariance matrix. (See Figure \ref{eigenvectormethod})
One way of looking at this approximation is that the difference 
of the sample spectrum to true spectrum can be viewed as a kind 
of moving average. This sort of moving average of the
sample spectrum involve some weights of the moving average which 
gradually modify as we change location. These weights are then given
by the columns of a matrix $B(\widehat{\lambda}_1,\widehat{\lambda}_2,\ldots,
\widehat{\lambda}_p)$, where $B(.)$ is the matrix defined in \eqref{Bij}.
From this point of view, our formulas have one very important advantage.  
Namely, they give a simple intuitive explanation of how approximately
the spectrum of the covariance gradually changes into the sample covariance. 
In contrast, the free probability approach is based on  algebraic properties 
which is convoluted enough and to our best knowledge alludes an intuitive 
explanation of what the transform ``does'' to the spectrum of the covariance. 
On the other hand, the free probability approach from a theoretical point
of view, has an advantage, namely, one can prove convergence
when $p$ goes to infinity when the spectrum of true covariance converges.
(But, that spectrum does not converge in measure since some eigenvalues
are of order $O(1)$ whilst others are $O(p)$ for real-life-data).
 Furthermore, a practitioner might
argue that since we have proven the spectrum of the true
covariance to be a fixed point of  our map \eqref{themap} 
(up to a smaller order term), even for finite $p$ we can get very good
approximations of the true spectrum.  We will study more mathematical properties of 
the map \eqref{themap} in an upcoming paper which will put on solid 
grounds this approach.

\subsection{Examples of $Z^tZ$ and their applications for machine learning}

Think of $Z$ as a matrix where each column represents a point in 
$\mathbb{R}^n$ and that we have a machine learning problem concerning
these points. Note that most machine learning algorithms depend only
on the relative position of the input data-points to each other.
If we know all the scalar products between the vectors representing
 the data-points, then we know the relative positions of the points
to each other. This is so because we can then determine the angles
between the vectors representing these points and also their length.
Thus, if the columns of $Z$ represent the data-points we want to use
for a machine learning problem, then  $Z^tZ$ gives all
the scalar products between the vectors representing
these points. Consequently, $Z^tZ$ defines in a unique way 
the relative position of the points to each other.
Most machine learning algorithms depend only on the relative position
of input-points to each other, and hence only on $Z^tZ$.
This shows why most machine learning algorithms can be thought of 
as having as input a matrix $Z^tZ$, where the columns of $Z$
represent the data-points for the machine learning task at hand.

Now,  one of the first things practitioners who look at 
a new data set do is to consider the spectrum of $Z^tZ$.
In reality however, they would be more interested in the spectrum
of $E[Z^tZ]$ and with very large feature space and not that large sample size,
the difference can be big. That is, the difference between
the spectrum of $Z^tZ$ and the spectrum
of $E[Z^tZ]$ is not negligible when $n$ is not much larger than $p$.

Let us give an example of a kernel for document classification,
in order to explain why the spectrum of the $Z^tZ$ is of interest to practitioners.

\begin{example}\label{e:ex1}
Consider that we have a collection of short documents
$D_1,D_2,\ldots,D_n$. Let $w_1,w_2,\ldots,w_p$
be a sequence of non-random words, which occur in
our documents. Take $Z$ to be the $n\times p$ matrix,
whose $ij$-th entry is $1$ if the $j$-th word occurs
in the $i$-th document and $0$ otherwise. 
The matrix $Z^tZ$ has then as $ij$-entry the 
number of documents where the words $w_i$ and $w_j$ co-occur.
Let $\vec{p}$ and $\vec{q}$ be two probability distributions
over the set of words. This simply means that $\vec{p}$ 
and $\vec{q}$ both have length $p$ and have non-negative entries 
which sum up to one. Assume that every documents gets randomly labeled
as topic 1 or topic 2 document. We can generate every document
as a sequence of i.i.d. words, drawn from either the topic 1 or topic
2 distribution depending on the topic classification
of the document. We describe here a simple way to think
about how the documents where generated by a random process.
Let $m$ be the number of words in each document.
Assuming $m$ is not too large, we get the approximation
\begin{equation}
\label{grossepute}
E[Z^tZ]\approx
n\cdot m\left(\;(m-1)(P(T_1)\cdot \vec{p}\otimes\vec{p}+P(T_2)\cdot\vec{q}\otimes\vec{q})+
diag(P(T_1)\cdot \vec{p}+P(T_2)\cdot \vec{q})\;\right)
\end{equation}

Where $diag(\vec{x})$ represents the diagonal $p\times p$-matrix
having on the diagonal the vector $\vec{x}\in \mathbb{R}^p$ and $T_1$, resp.
$T_2$ designates the first, respectively the second topic.

Note that the matrix  $\vec{p}\otimes\vec{p}$, respectively
$\vec{q}\otimes\vec{q}$
has as unique non-zero eigenvalue $|\vec{p}|^2$, resp $|\vec{q}|^2$.
Let $\vec{p}=(p_1,p_2,\ldots,p_p)$, then
\[
\max_i p_i\geq |\vec{p}|^2\geq \min_i p_i.
\]
and similarly
\[
\max_i q_i\geq |\vec{q}|^2\geq \min_i q_i.
\]

Thus, the matrix $E[Z^tZ]$ is up to linear factor a sum
of the two terms.  The first one 
\begin{equation}
\label{pT1}
(m-1)\left(P(T_1)\cdot \vec{p}\otimes\vec{p}+P(T_2)\cdot\vec{q}\otimes\vec{q}\right)
\end{equation}
and the second one is the diagonal matrix
\begin{equation}
\label{diag}
diag(P(T_1)\cdot \vec{p}+P(T_2)\cdot \vec{q}).
\end{equation}
As soon as $m$ is not too small, the eigenvalues of \eqref{pT1} will dominate 
the ones of \eqref{diag} due to the factor $(m-1)$.
Here, we have two topics only, but in a realistic model with 
$k$ topics and $m$ not too large we should see
$k$ eigenvalues of order $O(nm^2 |\vec{p}_j|^2)$
corresponding to the $k$ topics for $j=1,2,\ldots,k$ and the rest of the eigenvalues
will roughly correspond to the probabilities
in the probability vector for the words times $nm$.
Here $\vec{p}_j$ designated the probability distribution
of the words, given the $j$-th topic which we call $T_j$.
Thus, under this model, the spectrum
of the covariance matrix will contain $k$-higher order
eigenvalues corresponding to the topics and the rest
of the eigenvalues will be of order $O(1)$, so in 
the spectrum of the matrix $E[Z^tZ]$ we expect to
see an ``elbow'' at the $k$-th eigenvalue.
Ideally, for real data, looking at the spectrum should tell us how many
topics there are, more precisely this should be quantified by the ``elbow'' 
which tells us where the top eigenvalues start detaching from the rest of spectrum.
Now often you don't see an elbow in the spectrum
of $Z^tZ$ although there might be one in $E[Z^tZ]$, 
the reason being, that the spectrum of $Z^tZ$ can 
be considered a smoothing of the spectrum of $E[Z^tZ]$.
\end{example}

The previous example is  mainly here to show how the spectrum
of $E[Z^tZ]$ is important in most machine learning problems.
This is not the first situation we encountered.  
We started working on the reconstruction of the ground-truth
spectrum in the context of financial data that
is not too far from multidimensional normal distribution.
It is true that for short documents, the matrix $Z^tZ$
of the previous example is sparse. Our method, of spectrum
reconstruction does not yet work on sparse data.

With real data, often the spectrum
of the estimated covariance matrix, $Z^tZ$ has \emph{little fluctuation
in it, but  a big bias}, due to the concentration
of measure phenomena. Thus this spectrum can be considered almost as 
a non-random function of the spectrum of the true
covariance matrix. If this ``almost non-random map''
is injective it would be possible to  reconstruct almost 
exactly the true spectrum from the estimated one, at least assuming
that the data is multivariate normal.

Let us give another example. This is the example of stocks, which
is pedagogically a very good example although it is not where 
PCA is necessarily most useful.  Indeed, for stocks, the sectors 
and the capitalization are factors which are known and 
for which we do not need PCA to discover them. 

\begin{example}\label{e:ex2}
Let $\vec{Z}=(Z_1,Z_2,\ldots,Z_p)$
be a random vector with $0$ expectation
and multivariate normal probability distribution.
Let $Z$ be $n$ times $p$ matrix
where $Z_{ij}$ is the
daily return on day $i$ of the $j$-th stock
in a given portfolio. Hence we have $p$ stocks
and $n$ days considered. Let us denote by $\vec{Z}_i$
the $i$-th row of $Z$. We assume that
$\vec{Z},\vec{Z}_1,\vec{Z}_2,\ldots$ is a i.i.d.
sequence of random normal vectors with $0$ expectation.
The $0$ expectation comes from the fact that
 on the daily base the change in value of a stock
has an expectation of smaller order than its standard deviation.
Therefore in practice we can assume $0$ expectation. Because,
we assumed the rows to be i.i.d., this means that from one day
to the next there is no correlation, and there is a process
which is stationary in time. Hence the covariance of the 
$j$-th stock with the $k$-th stock on any given day,
is equal to
$$COV(Z_j,Z_k)=E[Z_j\cdot Z_k],$$
where we used the fact that $Z_j$ and
$Z_k$ have $0$ expectation. Consequently, the covariance
here is an expectation, and we estimate expectations usually
by taking average of independent copies of the variables
under consideration. The independent copies of $Z_j\cdot Z_k$
are given by taking these variables on the different days
$i=1,2,\ldots,n$. Thus, our estimate is the average
over different days
\begin{equation}
\label{COVZ}
\widehat{COV}(Z_j,Z_k)=\widehat{E}[Z_j\cdot Z_k]:=
\frac{Z_{1j}\cdot Z_{1k}+Z_{2j}\cdot Z_{2k}+\ldots+
Z_{nj}\cdot Z_{nk}}{n}
\end{equation}
for all $j,k\leq n$.  We are going to estimate 
every entry of the covariance matrix $COV[\vec{Z}]=(COV(Z_j,Z_k)_{j,k}$ using
formula \ref{COVZ}. This gives us the
following estimate of the covariance matrix
\begin{equation}
\label{estimatedcov}
\widehat{COV}[\vec{Z}]:=\frac{Z^t\cdot Z}{n}.
\end{equation}
We are going to consider a situation where both $n$ and $p$
go to infinity at the same time, whilst their ratio remains constant.
Now, when $p$ and $n$ are in a constant ratio, 
and go to infinity both at the same time, then
the estimate \eqref{COVZ} can be very bad. Indeed, if the sample 
size $n$ is half the size of the portfolio $p$, 
then, the matrix \eqref{estimatedcov} is defective and
half of the values in its spectrum are $0$.
This is certainly not the case for the true covariance matrix
$COV[\vec{X}]=(E[X_{1s}X_{1t}])_{s,t}$ since otherwise
there would be a lot of riskless investment opportunities.
Therefore, there is a substantial error, which means that 
there is a big difference between the spectrum of the original 
covariance matrix $\widehat{COV}[\vec{X}]$ and the estimated
one $\widehat{COV}[\vec{X}]$.
\end{example}

\subsection{The model and stability}

Let $\lambda_1>\lambda_2>\ldots>\lambda_p$ denote
the spectrum of the true covariance matrix $COV[\vec{Z}]$
and let $\widehat{\lambda}_1>\widehat{\lambda}_2>\ldots>\widehat{\lambda}_p$ denote
the spectrum of the estimated covariance matrix $\widehat{COV}[\vec{Z}]$.

\subsubsection{The first stability phenomena, under resimulation}
\begin{center}
\begin{figure}[!ht]
\caption{}
\label{comparespectrum500.jpeg}
\includegraphics[width=90mm]{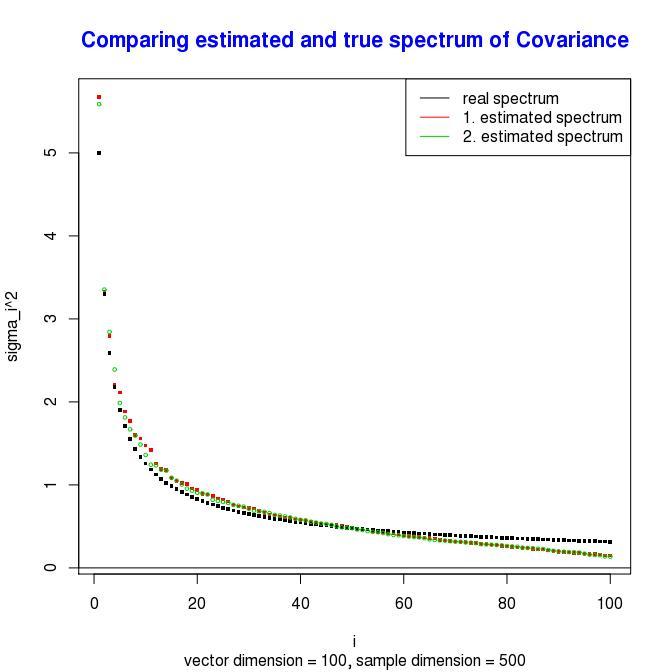}
\end{figure}
\end{center}

The main aim of this paper is to compare and explain methods
(including our own) for recovering the true spectrum $\lambda_i$, $i=1,2,\ldots,p$
from the estimates $\widehat{\lambda}_i$, $i=1,2,\ldots,p$.
Let us first continue with the framework of the Example~\ref{e:ex2} with a situation 
of some real data from stocks. This is shown in figure \ref{comparespectrum500.jpeg}
in which we took $p=100$ and $n=500$. In that picture, 
we have the true spectrum in black dots. Then, we resimulated
twice. Leading to the red line for the first simulation
and the green line for the second.  It is important to notice that 
there is a substantial difference between the black line and  the other 
two, the red and the green, however, the red and green curves are very close.  
This shows a situation where for practical purposes,
the estimated spectrum $\widehat{\sigma}_1,\widehat{\sigma}_2,\ldots,
\widehat{\sigma}_p$ has very little randomness in itself.

The spectrum of the estimated covariance matrix  
seems to be almost like a ``non-random function''
of the spectrum of the original covariance matrix.
One idea of recovering the true spectrum is the following. 
Take an arbitrary spectrum and simulate with this spectrum an estimated 
covariance matrix, then look at the spectrum and compare with the 
one we have at hand. Discard the original spectrum and then repeat this.  
However this seems very far from achieving any results at all.  We will proceed 
a little more methodical.  

\subsubsection{The second stability phenomena, under the rescaling 
or reduction/blowup}\label{s:stability}

\begin{center}
\begin{figure}[!ht]
\caption{}
\label{blowup1000}
\includegraphics[width=90mm]{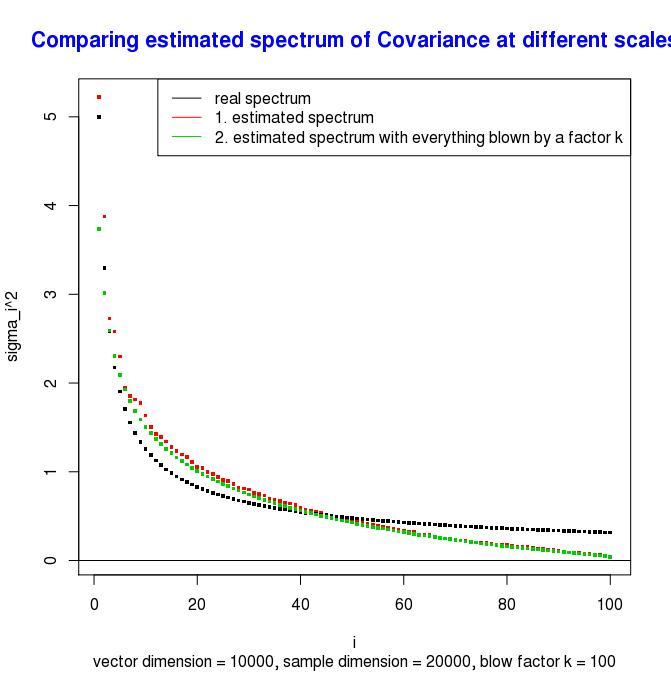}
\end{figure}
\end{center}

Another phenomena, we observed in the data
and our simulations is that the rescaled data's sample covariance spectrum 
does not change too much in distribution. Let us explain, assume we take
$n$ and $p$ to be a fixed proportion of a few times bigger such that their ratio
remains the same. Then, we take the same distribution for the spectrum
of the true covariance matrix for both simulations. 
In figure \ref{blowup1000}, we see the result of
such simulations. There, the red line is the sample covariance
spectrum when $n=200$ and $p=100$. The green line
is the sample covariance with everything $100$ times bigger.
That is, we take $n=20000$ and $p=10000$. The distribution
of the true spectrum is taken the same for both simulations
and is shown as the black line.  For the green line 
we took only one in $100$ of the eigenvalues for the plot.

Hence, from Figure \ref{blowup1000}, we see that the distribution
of the spectrum of the estimated covariance matrix
is almost unchanged by blowing up everything by a factor $100$
whilst keeping the same distribution for the spectrum
of the original covariance matrix. This tends to indicate that
we do not need $n$ to be very large to get close to the limiting
distribution.  Indeed the free-probability theory, as we explain in  Subsection
\ref{freeprobapproach},
predicts that if we keep the same distribution
for the spectrum of the original covariance matrix, then 
the spectrum of the estimated one converges as $n$ and $p$
go to infinity while we preserve their ratio.
However, the interesting thing we noticed with real data,
is that this phenomena of stability 
happens already with quite small data like here
in Figure \ref{blowup1000}, and not just 
with dimension size in the billions. 

Our spectrum reconstruction 
methods were first developed for when $Z$ has i.i.d rows 
with $0$ expectation which are multivariate normal. 
Our methods seems to work 
well even when the data is not exactly normal.
Now let $\lambda_1>\lambda_2>\ldots>\lambda_p$ denote
the eigenvalues of $E[Z^tZ]/n=COV[\vec{Z}]$, where again
$\vec{Z}$ designates a random vector having same distribution
then the rows of $Z$. Let $\widehat{\lambda}_1>\widehat{\lambda}_2>\ldots>
\widehat{\lambda}_p$ designate the eigenvalues of the sample covariance matrix
$Z^tZ/n$. The $i$-th estimation error is designated by $\Delta\lambda_i$
where $\Delta\lambda_i:=\widehat{\lambda}_i-\lambda_i$. In modern big data statistics
one is more interested in the case where $n$ and $p$ both go to infinity
at the same time, whilst their ratio $n/p$ converges to a constant.
The case of traditional statistics is when $p$ is fixed and $n$ goes to 
infinity. In that case, there is an asymptotic formula for 
$\lambda_i$ which is totally different from our case. 
In the finite dimension case, one takes  $\Delta\lambda$
to be asymptotically normal with expectation $0$ and a standard
deviation of order $O\left(\frac{1}{\sqrt{n}}\right)$.
It means that in this regime when $p$ is fixed and $n$ goes
to infinity alone, the expectation of $\Delta\lambda$ is negligible
compared to its standard deviation. In our case, the opposite holds, namely, 
the bias dominates the fluctuation so we can consider in a first approximation,
$\Delta\lambda$ to be non-random. We develop an approximation formula
for the eigenvectors and eigenvalues of the sample covariance.
We do this via the formula for covariance error matrix  \cite{Kolt1},
\cite{Kolt2}, \cite{Kolt3}
proven by Koltchinskii and Lounici, leading to the approximation formula 
\eqref{library} for the infinite
dimensional case ($n$ and $p$ both going to infinity
at the same time).

The second method is a fixed point method based on the eigenvectors
of both the covariance and the sample covariance matrix.

\subsubsection{Summary of the results}
\begin{itemize}
\item{}Our first spectrum reconstruction method is
based on the adaption of the finite dimensional
asymptotic formula for Principal Components to the infinite dimensional
case via the formula of Koltchinskii and Lounici \cite{Kolt1} for the spectral norm of covariance estimation error matrix.  The infinite dimensional
approximation formula is given by:
\begin{equation}
\label{library}
\Delta\lambda_i\approx\frac{\lambda_i}{n}\sum_{j\neq i} 
\frac{\widehat{\lambda}_j}{\widehat{\lambda}_j-\widehat{\lambda}_i},
\end{equation}
where $\Delta\lambda_i:=\widehat{\lambda}_i-\lambda_i$. Our first
method, now consists in ``solving''  for $\lambda_i$ 
given  all the sample covariance eigenvalues $\widehat{\lambda}_j$.  Thus we obtain 
\begin{equation}\label{approxlambda}
\lambda_{i}\approx \frac{\widehat{\lambda}_i}{1+\frac{1}{n}\sum_{j\neq i} 
\frac{\widehat{\lambda}_j}{\widehat{\lambda}_j-\widehat{\lambda}_i}}
\end{equation}

(For details on this method see subsection
\ref{aymptoticmethod}).
  \item{}Our second spectrum reconstruction method, consists
of finding  a fixed-point  for the map
$$\vec{\nu}\mapsto (\widehat{\lambda}_1,\widehat{\lambda}_2,\ldots,\widehat{\lambda}_p)
B(\vec{\nu}),$$
where $\vec{\nu}=(\nu_1,\nu_2,\ldots,\nu_p)$ with 
$\nu_1>\nu_2>\nu_3>\ldots>\nu_p\geq 0$ and $B(\vec{\nu})$ designating
the matrix with the expected coefficients square of the 
unit eigenvectors of the original covariance matrix expressed
in the basis of the unit eigenvectors of the sample covariance matrix.
This matrix $B(\vec{\nu})$ is for when the 
spectrum of true covariance is equal to $\vec{\nu}$.
(For details see Subsection \ref{SectionEigen}).
\end{itemize}

\subsubsection{Testing and consistency} 
\begin{center}
\begin{figure}[!ht]
\caption{Our reconstruction method for synthetic data}
\label{eigenvectormethod}
\includegraphics[width=90mm]{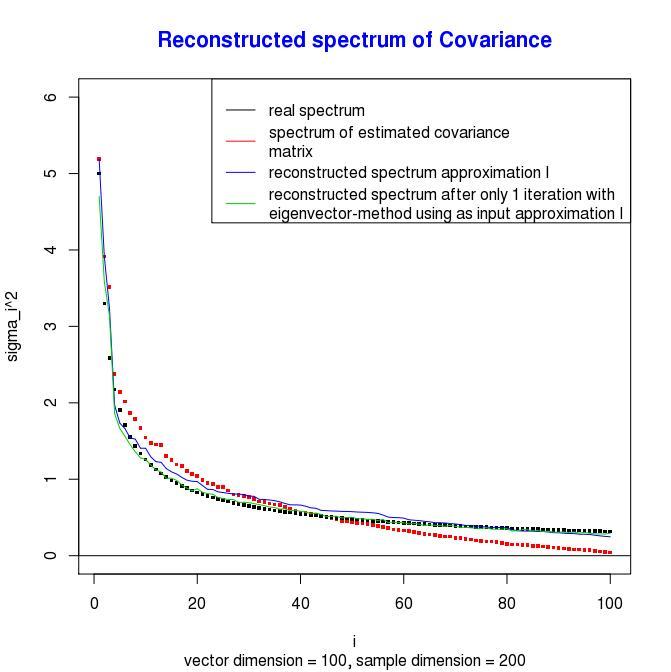}
\end{figure}
\end{center}

We first tested our spectrum reconstruction methods   
on synthetic multivariate normal data, with a real spectrum
taken from real life data and used simulations with this.
This means, that we simulate multivariate normal data using 
as ``more or less realistic'' spectrum  retrieved from real life data.
That is a spectrum $s\mapsto \lambda_s$, which is convex and has more or less
``the derivative going to $-\infty$ as $s$ approaches $0$''. 
Such kind of spectrum is what we noticed in real data.
This is what you can see in  Figure \ref{eigenvectormethod},
where we have taken the spectrum from real stocks and estimated it, by 
just taking the spectrum for the sample covariance when sample
size is much bigger than the vector size. In this way, we are pretty sure
that the sample spectrum is close to the true underlying spectrum.
This is so because when we hold $p$ fixed and let $n$ go to infinity,
then the sample spectrum converges to the true spectrum.

For this synthetic data with realistic spectrum, the result 
of our reconstruction is surprisingly good as soon as $n/p\ge2$
as can be seen in Figure \ref{eigenvectormethod}.

Below in Subsection \ref{testingdata}, we test our methods
and others on synthetic multivariate normal data and we see that except
when the spectrum is close to being flat (close to being constant), our 
fixed-point-eigenvector method
works really well. This only shows that when 
the data is truly hundred percent multivariate normal
our methods work well with many examples of spectrae. 
However, real data is never exactly multivariate normal, 
only approximately. Thus we need to test
our methods also on real data, and not just on synthetic data.

How can we be sure that it works with real data for recovering the ``true
underlying'' covariance matrix's spectrum, since we will never know
that spectrum exactly?  A first step, is to take 
a sample size much bigger than vector size 
(provided we have enough data for this).  
Then we see that increasing the sample size almost does not change the 
spectrum of the sample covariance matrix, we can assume that 
we have the true spectrum more or less.  Indeed for $p$ fixed and $n$ going 
to infinity the sample spectrum always converges to the true spectrum, 
even if the model is not multivariate
normal. 

\begin{center}
\begin{figure}[!h]
\caption{Our reconstruction methods with real data }
\label{realdata}
\includegraphics[width=100mm]{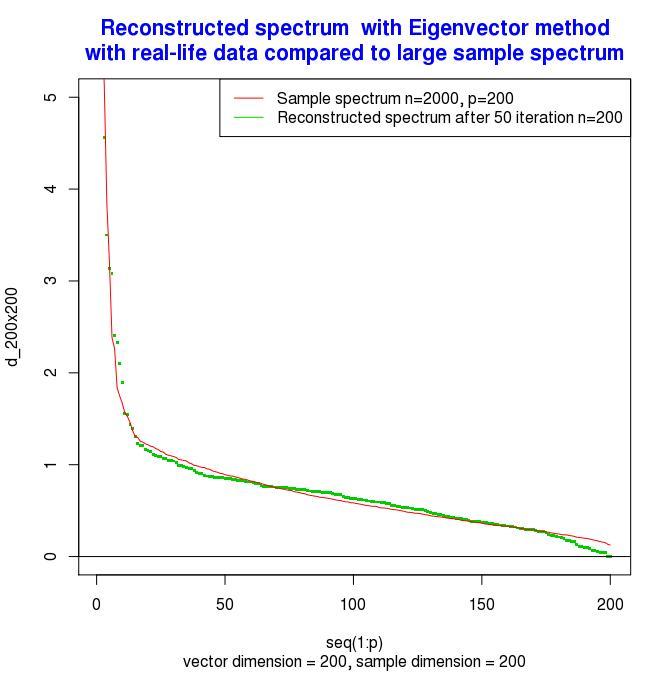}
\end{figure}
\end{center}

Thus, what we did is the following. We took two hundred stocks and their daily returns for 
$2000$ days. Then we took only $200$ days of real data (not synthetic data)
and applied our spectrum recovery method (eigenvalue-fixed point method)
to the real data with only 200 days.  We compare the reconstructed spectrum 
to the sample covariance matrix spectrum with $2000$
days, where we assume that the sample spectrum with $2000$ days
is close to the ``true underlying'' spectrum. The difference is quite
small indicating that our method does work in the 
multivariate normal data and also on real life  data. See for this 
Figure \ref{realdata}.

In real life, the data is never exactly i.i.d. For example with stocks
the true underlying covariance matrix will always change a little bit
over time. Thus, even with the perfect reconstruction method,
there must be a small discrepancy between the spectrum of the sample
covariance matrix for $n$ large and the reconstructed spectrum.
It is remarkable how little the difference between the two
is in Figure \ref{realdata}!

\begin{center}
\begin{figure}[!ht]
\caption{Reconstruction  with real data but two different
samples sizes compared }
\label{1000vs2000}
\includegraphics[width=90mm]{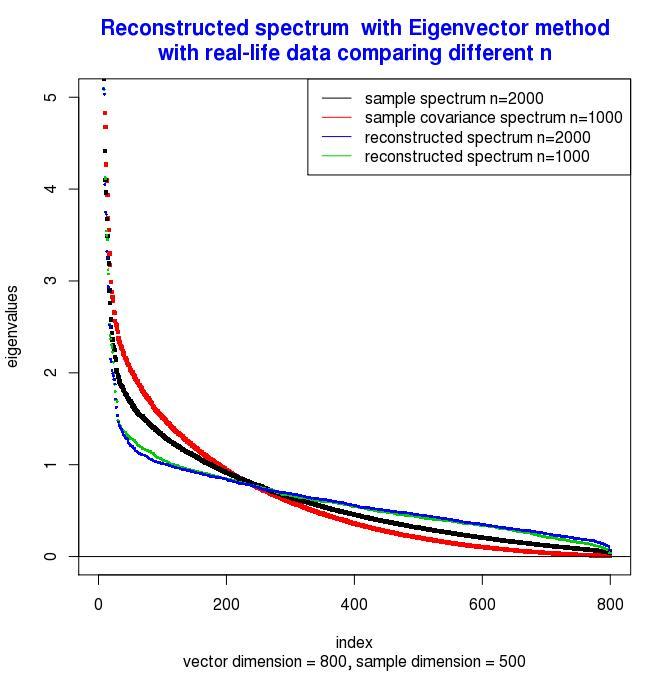}
\end{figure}
\end{center}

Another, way to find out if our reconstruction methods work
with real data, is to take a data set $Z$ of size $n$ times $p$.
Then, we do the reconstructions for different restrictions
of the sample size of our data set. That is we take 
$n_1<n$ and for the restricted data set where we take only $n_1$ samples
instead of $n$, that is for
$$(\left(Z_{i,j}\right)_{i\in[1,n_1],j}$$
we do the reconstruction and compare it to the reconstruction 
using all $n$ samples. This is what can be seen in Figure \ref{1000vs2000}.
There we took $800$ stocks and once the $2000$ days and another time only 
$1000$ days. We took $50$ iterations of the 
fixed-point-eigenvector method. The first input for the fixed
point method, was obtained by our asymptotic formula \ref{library}
for $\Delta\lambda$.

It is remarkable that though the spectrum of the sample
covariance matrix is very different between the two
cases $n=1000$ and $n=2000$, the reconstructed spectrum
are pretty close to each other This is a strong indication
that our reconstruction method works extremely well
for real data. Indeed, it seems very unlikely that for very different inputs 
(the two sample covariance spectrae
for $n=1000$ and the one for $n=2000$), would lead to almost same
output. The only plausible explanation for why the outputs are
about the same, is that they are close to the ``true underlying
spectrum'', which indeed is the same for both cases.

Note that to mitigate the fact that covariance matrix for
stocks changes over time, we did not take the first $1000$
days, when we reduced the sample size. Instead, we chose $1000$
days at random.

\subsection{Principal components}
\label{principalcomponents}
The eigenvalues of the covariance matrix $COV[\vec{Z}]$
are called {\it Principal Components} and we denote 
by $\vec{\mu}_1,\vec{\mu}_2,\ldots,\vec{\mu}_p$
 the eigenvectors (principal components)
of $COV[\vec{Z}]$, corresponding to the eigenvalues
$\lambda_1,\lambda_2,\ldots,\lambda_p$.
We assume that the eigenvalues are ordered as follows 
\[
\lambda_1>\lambda_2>\ldots>\lambda_p. 
\]

The quantities associated to the estimated covariance matrix are distinguished by a hat.  
Thus, the eigenvectors of the 
estimated covariance matrix $\widehat{COV}[\vec{Z}]$
are denoted by
\[
\widehat{\vec{\mu}}_1,\widehat{\vec{\mu}}_2,\ldots,\widehat{\vec{\mu}}_p,
\]
with corresponding eigenvalues:
\[
\widehat{\lambda}_1,\Hat{\lambda}_2,\ldots,\widehat{\lambda}_p.
\]
Why does one use eigenvectors? At least with some theoretical models
the principal components with large eigenvalues correspond to
factors. Let us see an example from finance which is related to our already mentioned 
Example \ref{e:ex2}.  This model is not necessarily
the one where principal components are most useful, but
is in our opinion one of the best for pedagogical reasons.

{\footnotesize One good example to understand where our $3$-order model
for eigenvalues come from, are financial
 stocks. Assume that 
$\vec{Z}=({\color{red}X},{\color{red}Y},{\color{green}Z},{\color{green}W})$ denotes
the daily change of value of four stocks.
\begin{itemize}
\item{}Stocks ${\color{red}X}$ and ${\color{red}Y}$ depend only on index 
${\color{red}T}$
and stocks ${\color{green}Z}$ and ${\color{green}W}$ depend only on index 
${\color{green}S}$.
Thus, 
\begin{align*}
&{\color{red}X}=a_Z{\color{red}T}+0\cdot S+\epsilon_X\\
&{\color{red}Y}=a_Y{\color{red}T}+0\cdot S+\epsilon_Y\\
&{\color{green}Z}=0\cdot T+b_Z{\color{green}S}+\epsilon_Z\\
&{\color{green}W}=0\cdot T+b_W{\color{green}S}+\epsilon_W
\end{align*}

\item{}We assume that ${\color{green}S}$, 
${\color{red}T}$ and $\epsilon_X,\epsilon_Y,\epsilon_Z,\epsilon_W$
are all independent of each other. The firm specific terms
$\epsilon_X,\epsilon_Y,\epsilon_Z,\epsilon_W$ are supposed to all
have the same  st. deviation $\sigma$. 
Furthermore $a_Z,a_Y,b_Z,b_W$ are non-random
and let $VAR[S]=VAR[T]=1$\\
\end{itemize}

With this finite sector model we get:
$$COV[\vec{Z}]=
\left(
\begin{array}{cccc}
{\color{red}a_X^2}&{\color{red}a_Za_Y}&0&0\\
{\color{red}a_Ya_Z}&{\color{red}a_y^2}&0&0\\
0&0&{\color{green}b_Z^2}&{\color{green}b_Zb_W}\\
0&0&{\color{green}b_Wb_Z}&{\color{green}b_w^2}
\end{array}
\right)+\sigma^2 I
$$
where $I$ is identity matrix.
The  two eigenvectors with  biggest eigenvalues are
$${\color{red}\vec{\mu}_1=(a_Z,a_Y,0,0)}\;\;,\;\;
{\color{green}\vec{\mu}_2=(0,0,b_Z,b_W)}$$ with eigenvalues
${\color{red}a_Z^2+a_y^2+\sigma^2}$, respectively 
${\color{green}b_Z^2+b_W^2+\sigma^2}$.
We can easily extend this model to higher dimensional case.  
Assume for this that we have $2p$ stocks. The first $p$ 
are in the first sector and the second half 
in the second sector. Let 
$$\vec{Z}=({\color{red}Z_1},{\color{red}Z_2},\
\ldots,{\color{red}Z_p},{\color{green}Z_{p+1}},
{\color{green}Z_{p+2}},\ldots,{\color{green}Z_{2p}})$$ denote
the daily change of value of our $2p$ stocks summarized in one vector.
With this two sector model 
$$Z_i=a_i {\color{red}T}+\epsilon_i$$ for 
$i={\color{red}1,\ldots,p}$
and $$Z_j=b_j {\color{green}S}+\epsilon_j$$ for $j={\color{green}p+1,\ldots,2p}$ We assume
 ${\color{green}S}$,${\color{red}T}$ and $\epsilon_i$'s independent
of each other with $0$ expectation. At first assume all $\epsilon_i$ have
same standard  deviation $\sigma$.
We find then the following covariance matrix
\[
COV[\vec{Z}]=
\left(
\begin{array}{ccccccc}
{\color{red}a_1^2}&{\color{red}a_1a_2}&\ldots&{\color{red}a_1a_p}&0&\ldots&0\\
{\color{red}a_2a_1}&{\color{red}a_2^2}&\ldots&{\color{red}a_2a_p}&0&\ldots&0\\
&&&\ldots&&&\\
&&&\ldots&&&\\
{\color{red}a_pa_1}&{\color{red}a_pa_2}&\ldots&{\color{red}a_pa_p}&0&\ldots&0\\
0&0&\ldots&0&{\color{green}b_{p+1}b_{p+1}}&\ldots&{\color{green}b_{p+1}b_{2p}}\\
0&0&\ldots&0&{\color{green}b_{p+2}b_{p+1}}&\ldots&{\color{green}b_{p+2}b_{2p}}\\
&&&\ldots&&&\\
&&&\ldots&&&\\
0&0&\ldots&0&{\color{green}b_{2p}b_{p+1}}&\ldots&{\color{green}b_{2p}b_{2p}}\\
\end{array}
\right)+\sigma^2 I
\]

With the model just introduced, the  two eigenvectors with  biggest eigenvalues are
\[
{\color{red}\vec{\mu}_1=(a_1,a_2,\ldots,a_p,0,\ldots,0)}\;\;,\;\;
{\color{green}\vec{\mu}_2=(0,\ldots,0,b_{p+1},b_{p+2},\ldots,b_{2p})}
\]
with eigenvalues
\[
{\color{red}a_1^2+a_2^2+\ldots+a_p^2+\sigma^2}=O(p),
\]
respectively 
\[
{\color{green}b_{p+1}^2+b_{p+2}^2+\ldots+b_{2p}^2+\sigma^2}=O(p).
\]
 
Other eigenvalue is $\sigma=O(1)$ with multiplicity $p-2$.
Typically we take the data standardized (i.e. a correlation matrix)
and the stocks have standard deviation $1$. This means
that the coefficients $a_i$ and $b_i$ are correlation coefficients
and are hence of order $O(1)$.  In real life we have 
the sector factors (like here in the example $S$ and $T$) are not independent
but depend all on the general economy. Hence, more realistically,
we have a linear system of equation with factors, but the factors
are not independent, then the vectors with the {\it loadings} of the factors
$${\color{red}(a_1,a_2,\ldots,a_p,0,\ldots,0)}\;\;,\;\;
{\color{green}(0,\ldots,0,b_{p+1},b_{p+2},\ldots,b_{2p})}$$
are no longer the principal components. However, we still 
have that the span of these two factor-vectors is the same
as the span of the two first principal components $\vec{\mu}_1$
and $\vec{\mu}_2$. (At least if we assume that
the $\epsilon_j$'s all have same variance and are independent
of the factors $S$ and $T$ and of each other).
Thus, in that case, the principal components
are not interesting separately, but only in group.
This means we are interested in their span.}

If we can retrieve the span of the eigenvectors
with large eigenvalues, we will
use that span for dimension reduction and then
look for the factors within that span. This is called 
{\it   rotation of factors} which essentially tries to retrieve
the factors from their span. Nevertheless, we see
in our example that typically, the eigenvalues corresponding to
factors are of order $O(p)$. The other eigenvalues
coming from the terms $\epsilon_i$, are of order $O(1)$.
Thus, if we were given the ``true'' spectrum 
$\lambda_1,\lambda_2,\ldots,\lambda_p$, we would try to see 
where the big eigenvalues ``end'', and where
the smaller eigenvalues ``start'' in the spectrum. In
reality this is often very difficult to see from the estimated
spectrum $\widehat{\lambda}_1,\widehat{\lambda}_2,\ldots,\widehat{\lambda}_p$
where the ``big eigenvalues stop'' and where the ``small eigenvalues
start''.  This is so because if $n$ is not much larger than $p$
there is a big estimation error in the spectrum of the covariance
matrix. This is to say that there can be a big difference between
\begin{equation}
\label{truespec}
\lambda_1,\lambda_2,\ldots,\lambda_p,
\end{equation}
and 
\begin{equation}
\label{estspec}
\widehat{\lambda}_1,\widehat{\lambda}_2,\ldots,\widehat{\lambda}_p,
\end{equation}
when $n$ is not much bigger than $p$. (In our practical experience
we need $n$ to be about $10$ times bigger than $p$ in order
to be able to consider that the estimated spectrum
\ref{estspec} can be considered identical to
the true spectrum \ref{truespec} in practice).
Therefore the usual error when taking the estimated spectrum
\ref{estspec} instead of the true spectrum \ref{truespec}
is a sort of smoothing. (We will describe below why)
Due to this, very often, in the original spectrum 
we can really see where the ``large eigenvalues end in the 
true spectrum'' which often is no longer easily visible in the estimated
spectrum \ref{estspec} due to the smoothing. Indeed, the 
passage corresponds to an ``elbow'' in the original
spectrum \ref{truespec} which is no longer visible
in the estimated spectrum \ref{estspec}. 

That is one should be able to better identify
where the eigenvalues corresponding to factors end and where
in the spectrum, the noise starts. Another motivation, is that
with normal data, if we are given the original spectrum \ref{truespec},
we can simulate everything. This can then allow as we explain
below to see if the estimated eigenvectors are typically
very different from the true ones.

\section{Our recovering methods and the free probability
approach
explained}

We explain next our two  algorithms and where they come from.
Given as input the spectrum of $\widehat{COV}[\vec{X}]$
$$\widehat{\lambda}_1,\widehat{\lambda}_2,\ldots,\widehat{\lambda}_p$$
we try to recover the  spectrum of
the original covariance matrix $COV[\vec{X}]$
\[
\lambda_1,\lambda_2,\ldots,\lambda_p.
\]
In the last subsection of this section we explain the general
idea behind the free probability approach.

Let us summarize the different methods.
\begin{itemize}
\item{}Free probability based approach. This approach has first been
 pioneered by  EI Karoui \cite{elkaroui2008},
was further developed by Silverstein \cite{silverstein1995}, Bai and Yin \cite{bai1988}, Yin, Bai and Krishnaiah \cite{yin1983}
and recently Ledoit and Wolf  \cite{ledoit2012}, \cite{li2013} and \cite{ledoit2015}.
\item{} Moment methods. \cite{Valiant}. This is the same methods as is usually
for proving the semi-circular law. But, then
for the higher moments, too much computation is needed, but we definitely
would use moments which are not too high,
when dealing with big real life data as a double insurance.

\item{}Modeling ``true'' spectrum  with low
dimension parameter function family and then apply
asymptotic statistics.

\item{} Shrinkage method which was pioneered by Stein \cite{Stein}, 
see also Bickel and Levina \cite{Bickel} or Donoho \cite{Donoho}.

\item{} Our methods.
 \begin{enumerate}

 \item{} Finite dimensional asymptotic formula
adapted to the infinite dimensional case using the formula
of Koltchinskii and Lounici \cite{Kolt1} for covariance error matrix spectral norm.
 
  \item{} Fixed point of map defined with the matrix of second moments
of principal component eigenvectors. For this we take
the eigenvectors of the true covariance matrix expressed
in the basis of the sample covariance eigenvectors as a function
of a given spectrum $\vec{\nu}=(\nu_1,\nu_2,\ldots,\nu_p)$.
      \end{enumerate}
\end{itemize} 
Let us next present our two  methods more in detail:

\subsection{Our first  method based on asymptotic finite dimensional
formula adapted via Koltchinski and Lounici}
\label{aymptoticmethod} 

  This method with real non-sparse data works quite well as long as 
$p/n$ is not much bigger than $1$. It is a simple formula,
which gives an approximation to
the correcting term to obtain the values $\lambda_i$
from the eigenvalues $\widehat{\lambda}_i$, $i=1,2,\ldots,p$
 of the estimated covariance matrix. The formula for approximating the
$i$-th error between the eigenvalue of the covariance matrix, and the corresponding
eigenvalue of the estimated covariance matrix is given. This error
is $\Delta\lambda_i:=\widehat{\lambda}_i-\lambda_i$ can be approximated as
in \ref{library}. One can then ``solve'' the approximation \eqref{library}
for $\lambda_i$ to obtain the approximation for the true eigenvalue $\lambda_i$
given in \eqref{approxlambda}. But, where does the formula \eqref{library} come from
in the first place? In an upcoming paper, we will present exact calculations,
which are rather elaborated. Here we can already show the main idea
and how the result of Koltschinskii and Lounici \cite{Kolt1} is essential
in getting the formula. We start with a three dimensional example:\\[3mm]

{\footnotesize
Assume we have a  three dimensional random vector with independent
normal entries $\vec{X}=(X,Y,Z)$ having each expectation $0$.
 We assume that
the standard deviations are decreasing: $\sigma_X>\sigma_Y>\sigma_Z$.
The covariance matrix we consider is:
\begin{equation}
\label{covvecX}COV[\vec{X}]=\left(
\begin{array}{ccc}
\sigma_X^2&0&0\\
0&\sigma_Y^2&0\\
0&0&\sigma_Z^2
\end{array}
\right).
\end{equation}
Now, for that covariance matrix, the eigenvectors are the three canonical
vectors: $(1,0,0)$, $(0,1,0)$ resp. $(0,0,1)$ with eigenvalues
$\sigma_X^2$, $\sigma_Y^2$ and $\sigma_Z^2$. We are going to look
at what happens to the eigenvectors, when instead of the covariance matrix
we take the estimated covariance matrix. Let us explain more in detail the
setting. 

Let us assume that $\vec{X},\vec{X}_1,\vec{X}_2,\ldots$ are independent copies
of the vector $\vec{X}$, where $\vec{X}_i=(X_i,Y_i,Z_i)$. We assume
that we have $n$ such vectors and want to estimated the covariance matrix 
\eqref{covvecX}. Let $X$ be the $n\times 3$ matrix, where the 
$i$the row of $X$ is equal to $\vec{X}_i$.
The estimated covariance matrix
is then given by
$$\widehat{COV}[\vec{X}]:=\frac{X^t\cdot X}{n}$$
The difference between true and estimated covariance matrix is denoted by $E$
so that
$$E=COV[\vec{X}]-\frac{X^t\cdot X}{n}.$$
Instead of the true covariance matrix $COV[\vec{X}]$
we consider the perturbated matrix
$$\widehat{COV}[\vec{X}]=COV[\vec{X}]+E$$
which is obtained from the true covariance matrix by adding $E$.

There is a general equation for what happens with eigenvectors
 when we add a perturbation $E$ to a matrix $A$. Say that $\vec{\mu}$ 
is an unitary eigenvector of $A$ with eignevalue $\lambda$. 
Let $\alpha \vec{\mu}+\Delta\vec{\mu}$
be a unitary eigenvector of $A+E$ with eigenvalue $\lambda+\Delta\lambda$.
Then, we  have
\begin{equation}
\label{I}
A\vec{\mu}=\lambda\vec{\mu}
\end{equation}
and also
\begin{equation}
\label{II}(A+E)(\alpha\vec{\mu}+\Delta\vec{\mu})=(\lambda+\Delta\lambda)(\alpha\vec{\mu}+\Delta\vec{\mu}).
\end{equation}
Subtracting $\alpha$ times \eqref{I} from \eqref{II} we find
\begin{equation}
\label{III}
(A-(\lambda+\Delta\lambda)I+E)
\Delta\vec{\mu}=-\alpha E\cdot \vec{\mu}+\alpha\Delta\lambda \cdot \vec{\mu}
\end{equation}

Now we can rewrite equation \eqref{III}  with $A$ being the covariance matrix
\ref{covvecX}. We consider the eigenvector $\vec{\mu}=(1,0,0)$
of $A$ with eigenvalue $\lambda=\sigma_X^2$. We take 
$\Delta\mu$ perpendicuar to $\vec{\mu}$. In the present,
case this means that $\Delta\mu=(0,\Delta\mu_Y,\Delta\mu_Z)$.
Thus, without leaving out terms, we find the following  exact equation:

\begin{align*}
\label{zaza3}
&\left(
\begin{array}{ccc}
0&0&0\\
0&\sigma_Y^2-\sigma_X^2-\Delta\lambda&0\\
0&0&\sigma_Z^2-\sigma_X^2-\Delta\lambda
\end{array}
\right)
\left(\begin{array}{c}
0\\
\Delta\mu_Y\\
\Delta\mu_Z
\end{array}\right)
=\\
&=
-\alpha\left(
\begin{array}{c}
E_{11}\\
E_{21}\\
E_{31}
\end{array}
\right)+
\alpha\left(
\begin{array}{c}
\Delta\lambda\\
0\\
0\\
\end{array}
\right)-
\left(
\begin{array}{ccc}
0&0&0\\
0&E_{22}&E_{23}\\
0&E_{32}&E_{33}
\end{array}\right)
\left(
\begin{array}{c}
0\\
\Delta\mu_Y\\
\Delta\mu_Z
\end{array}
\right)-\left(
\begin{array}{ccc}
0&E_{12}&E_{13}\\
0&0&0\\
0&0&0
\end{array}\right)
\left(
\begin{array}{c}
0\\
\Delta\mu_Y\\
\Delta\mu_Z
\end{array}
\right)
\end{align*} 
the above equation for matrices can be ``separated into two parts''.
First the single equation for $\Delta\lambda$:
\begin{equation}
\label{Deltalambda}
\Delta\lambda=E_{11}
+\frac{1}{\alpha}(E_{12}\Delta\mu_Y+
E_{13}\Delta\mu_Z).\end{equation}

Then the $p-1$ dimensional equation for $\Delta\vec{\mu}$ given
as follows:
\begin{align*}&\left(\left(
\begin{array}{cc}
\sigma_Y^2-\sigma_X^2-\Delta\lambda&0\\
0&\sigma_Z^2-\sigma_X^2-\Delta\lambda
\end{array}
\right)+\left(
\begin{array}{cc}
E_{22}&E_{23}\\
E_{32}&E_{33}
\end{array}\right)\right)
\left(\begin{array}{c}
\Delta\mu_Y\\
\Delta\mu_Z
\end{array}\right)
=\\
&=-\alpha
\left(
\begin{array}{c}
E_{21}\\
E_{31}
\end{array}
\right)
\end{align*} 

If $\Delta\lambda$ is given we can solve the above equation
for $\Delta\vec{\mu}=(\Delta\mu_Y,\Delta\mu_Z)$.
We find
$$\left(
\begin{array}{c}
\Delta \mu_Y\\
\Delta\mu_X
\end{array}
\right)=-\alpha\cdot
(I+D_1\cdot E_1)^{-1}\cdot D_1 \vec{E}_1=-\alpha(D_1^{-1}+E_1)^{-1}\vec{E}_1,$$
where
$$D_1=\left(
\begin{array}{cc}
\sigma_Y^2-\sigma_X^2-\Delta\lambda&0\\
0&\sigma_Z^2-\sigma_X^2-\Delta\lambda
\end{array}
\right)^{-1}$$
and $E_1$ is the reduce matrix obtained from $E$ by deleting the first row and first column:
$$E_1=\left(
\begin{array}{cc}
E_{22}&E_{23}\\
E_{32}&E_{33}
\end{array}\right).$$
Also, $\vec{E}_1$ is the first column of the matrix $E$ where we delete the first entry.
$$\vec{E}_1=\left(
\begin{array}{c}
E_{21}\\
E_{31}
\end{array}
\right).$$
Now, in general with a higher dimensional situation and say you want to estimate the $i$-th eigenvector.
You get the same type formula:
\begin{equation}
\label{YOU}
\Delta\vec\mu_i=-\alpha
(I+D_i\cdot E_i)^{-1}\cdot D_i \vec{E}_i,
\end{equation}
where $E_i$ is the reduced matrix obtained from $E$ by deleting the $i$ row and column.
Furthermore, $D_i$ is the diagonal $p\times p$-matrix, obtained by deleting the $i$-th row
and column from the diagonal matrix having as $j$ entry in the diagonal $1/(\sigma_j^2-\sigma_i^2-\delta\lambda_i)$.
Again $\vec{E}_i$ is the vector obtained from taking the $i$-th column of $E$ and deleting the $i$-th entry.
Now, we assume that the spectral norm of 
$D_iE_i$ is much less than one, so that
$${\tt assumption:} \;\;|D_iE_i|<<1$$
In that case, we can leave out in formula \ref{YOU}
the matrix $D_iE_i$ to obtain the approximation
$$\Delta\vec\mu_i\approx-\alpha
 D_i \vec{E}_i.$$
We plug the last approximation back into
\ref{Deltalambda} to obtain: 
\begin{equation}
\label{dumbo}
\Delta\lambda\approx
-\left(\frac{E_{12}^2}{\sigma^2_Y-\sigma^2_X-\Delta\lambda}+
\frac{E_{13}^2}{\sigma^2_Z-\sigma^2_X-\Delta\lambda}\right),
\end{equation}
where we also left out $E_{11}$ since it is a smaller term
of order $O(\frac{\sigma_X^2}{\sqrt{n}})$.
Let 
$$N_{12}=\frac{E_{12}}{\sigma_X\sigma_Y}\cdot\sqrt{n},
N_{13}=\frac{E_{13}}{\sigma_X\sigma_z}\cdot\sqrt{n}.$$
By central limit theorem, $N_{12}$ and $N_{13}$ are close
to standard normal. On top of it they are uncorrelated.
Hence, we can write equation \ref{dumbo} as
$$
\Delta\lambda\approx
-\frac{\sigma_X^2}{n}
\left[(\frac{N_{12}^2\sigma_Y^2}{\sigma^2_Y-\sigma^2_X-\Delta\lambda}+
\frac{N_{13}^2\sigma_Z^2}{\sigma^2_Z-\sigma^2_X-\Delta\lambda}\right],
$$
In the real data sets we have in mind, typically $\Delta\lambda$'s
standard deviation is of smaller order than $\Delta\lambda$.
This allows to take the expectation on the right side
of the last approximation above to find:
$$
\Delta\lambda\approx
-\frac{\sigma_X^2}{n}
\left[\frac{\sigma_Y^2}{\sigma^2_Y-\sigma^2_X-\Delta\lambda}+
\frac{\sigma_Z^2}{\sigma^2_Z-\sigma^2_X-\Delta\lambda}\right],
$$
where we acted as if $\Delta\lambda$ would be a constant and we used
the fact that $E[N_{12}^2]=E[N_{13}^2]=1$. The values $\sigma_Y^2$
and $\sigma_Z^2$ are not known, so we replace them by their estimates
leading to our formula:
\begin{equation}
\label{3D}
\Delta\lambda\approx-
\frac{\sigma_X^2}{n}
\left[\frac{\widehat{\sigma}_Y^2}{\widehat{\sigma}^2_Y-\widehat{\sigma}^2_X}+
\frac{\widehat{\sigma}_Z^2}{\widehat{\sigma}^2_Z-\widehat{\sigma}^2_X}\right],
\end{equation}
where we also used that $\widehat{\sigma}_X^2=\sigma_X^2+\Delta\lambda$.

With a bigger normal random vector
$$\vec{X}=(X_1,X_2,\ldots,X_p)$$
with independent normal entries with $0$ expectation
and 
$$\lambda_i=VAR[X_i]=\sigma_i^2$$
for all $i=1,2,\ldots,p$ we can generate a $p$-dimensional
random vector and then we can derive similarly to \eqref{3D} 
as in the three dimensional case
\begin{equation}
\label{pD}
\Delta\lambda_i:=\widehat{\lambda}_i-\lambda_i=
\widehat{\sigma}_i^2-\sigma_i^2\approx-
\frac{\sigma_i^2}{n}\sum_{j\neq i}\frac{\widehat{\sigma}_j^2}{\widehat{\sigma}_j^2-\widehat{\sigma}_i^2}.
\end{equation}
But this is just our approximation formula given in \eqref{library}.
The fact that we took a normal vector with independent entries
is not a restriction.  When we take as basis the principal component
of the covariance matrix, the components become independent,
so that formula \eqref{pD} holds not just for
normal with independent entries but also for normal vector with
entries not independent of each other.

The only problem for proving this formula in $p$-dimensional
case, is bounding the spectral norm of  the matrix product $D_iE_i$.
For this we note that this spectral norm is equal to the spectral
norm of 
\begin{equation}
\label{DED}|D_i|^{1/2}E_i|D_i|^{1/2}
\end{equation}
where $|D_i|^{1/2}$ represents the diagonal matrix obtained from $D_i$
by replacing each diagonal entry by the square root of the absolute value
for the corresponding entry.  The matrix given in \eqref{DED}
can be viewed as a error-sample covariance matrix,
but when the ``true'' spectrum is modified.
The modification consists of taking 
$\sigma^2_j/(\sigma^2_j-\sigma_j^2)$ instead of $\sigma_j^2$
for all $j$ eigenvalue of a true covariance matrix, for all
$j\neq i$. We can bound the spectral norm of \eqref{DED} with the result
of Koltchinskii and Lounici \cite{Kolt1} for covariance 
error matrix spectral norm. This is why the result
of Koltchinskii and Lounici is so
important to us.  As we already mentioned, the detail of these calculations will be explained more in detail
 in an upcoming paper, however the main idea is just explained above.  
}

\subsection{Our second  method: fixed-point for eigenvector-entry-square map}
\label{SectionEigen}
Take a vector of non-random coefficients
$$\vec{a}=(a_1,a_2,\ldots,a_p).$$
Consider the scalar product
\begin{equation}
\label{scalarproduct}
\vec{a}\cdot \vec{Z}=a_1Z_2+\ldots+a_pZ_p.
\end{equation}
The scalar product \eqref{scalarproduct} has expectation $0$.
To estimate its variance since the expectation
is $0$, we can take the average of independent copies square
\begin{equation}
\label{VARa}\widehat{VAR}[\vec{a}\cdot \vec{Z}]:=\widehat{E}[\left(\vec{a}\cdot
\vec{Z}\right)^2]=
\frac{(\vec{a}\cdot \vec{Z}_1)^2+\ldots+(\vec{a}\cdot \vec{Z}_n)^2}{n}
=\vec{a}^t\frac{Z^tZ}{n} \vec{a},
\end{equation}
where  $\vec{Z}_i$ denotes the $i$-th row of the $n\times p$-matrix
$Z$. Recall that the $\vec{Z}_1,\vec{Z}_2,\ldots$ is an i.i.d.
sequence of independent copies of $\vec{Z}$.
Thus our estimate is a mean of $n$ nicely behaved (
exponential decaying tail) independent variables.
Hence, the estimation error of \eqref{VARa} is with high probability of order 
$O(\frac{1}{\sqrt{n}})$. This is to say that
\begin{equation}
\label{en}
\frac{(\vec{a}\cdot \vec{Z}_1)^2+\ldots+(\vec{a}\cdot \vec{Z}_n)^2}{n}
=E[(\vec{a}\cdot \vec{Z}_1)^2]+O\left(\frac{\sigma_j^2}{\sqrt{n}}\right).
\end{equation}
Now, take $\vec{Z}$ to be a row verctor
and $\vec{a}$ to be a column vector. Then
\begin{equation}
\label{aZa}
E[(\vec{a}\cdot \vec{Z})^2]=
E[\vec{a}^t\vec{Z}^t\vec{Z}\vec{a}]=\vec{a}^t E[\vec{Z}^t\vec{Z}]\vec{a}=
\vec{a}^t COV[\vec{Z}]\vec{a}
\end{equation}
and so if $\vec{a}$ is the $j$-th unit-eigenvector of the original
covariance matrix, with eigenvalue $\lambda_j$, then the  expression on
the right most side of equation \ref{aZa} is equal to
\begin{equation}
\label{aZa2}
\vec{a}^t COV[\vec{Z}]\vec{a}=\lambda_j.
\end{equation}
In that case combining this with \eqref{aZa}, yields
$$
E[(\vec{a}\cdot \vec{Z})^2]=\lambda_j$$
and hence combining this with  \ref{VARa} and \ref{en}
\begin{equation}
\label{az}
\vec{a}^t\frac{Z^tZ}{n} \vec{a}=\lambda_j+O\left(\frac{\sigma_j^2}{\sqrt{n}}\right).
\end{equation}
The last equation above tells us that if we would know the $j$-th unit
eigenvector of the original covariance matrix $COV[\vec{Z}]$,
then we would be able to reconstruct the $j$-th eigenvalue $\lambda_j$
up to a smaller order term. Now let $U$ be the unitary matrix
whose columns are the unitary eigenvectors of the matrix $\frac{Z^tZ}{n}$
in the order of decreasing eigenvalues. Then, we can rewrite \eqref{az} as follows
\[
\lambda_j=\vec{a}^tU(U^t\frac{Z^tZ}{n}U)U^t\vec{a}+O\left(\frac{\sigma^2_j}{\sqrt{n}}\right).
\]
Now, $U^t\frac{Z^tZ}{n}U$ is the diagonal matrix of estimated eigenvalues $\widehat{\lambda}_j$.
(Recall that $\widehat{\lambda}_j$ designated the $j$-th eigenvalue of the estimated covariance
matrix $\widehat{COV}(\vec{Z})=\frac{Z^tZ}{n}$.) Furthermore $U^t\vec{a}$
is the $j$-th eigenvector of the true covariance matrix $COV[\vec{Z}]$ espressed in the basis of the eigenvectors
of $\frac{Z^tZ}{n}$.
So we have
\begin{equation}
\label{lambda}
\lambda_j= \vec{a}^tU  
\left(\begin{array}{cccccc}
\widehat{\lambda}_1&0&0                            & \ldots&0\\
0              &\widehat{\lambda}_2&0              &  \ldots&0\\
0              & 0             &\widehat{\lambda}_3&  \ldots &0\\
\ldots\\
0               &0              &0              &\ldots   &\widehat{\lambda}_p  
\end{array}
\right)U^t\vec{a}+O\left(\frac{\sigma^2_j}{\sqrt{n}}\right).
\end{equation}

Consequently, the estimated eigenvalues $\widehat{\lambda}_j$
have a fluctuation of smaller order than their own values and thus 
it makes sense to assume that 
\[
O(\sigma_{\widehat{\lambda}_j})<<O(E[{\widehat{\lambda}_j}]).
\]

What this means is that we can replace in \eqref{lambda} the 
$\widehat{\lambda}_i$'s by their expectation. 
This gives an error of smaller order. Hence we find that
\begin{equation}
\label{lambda2}
\lambda_j= \vec{a}^tU  
\left(\begin{array}{cccccc}
E[\widehat{\lambda}_1]&0&0                            & \ldots&0\\
0              &E[\widehat{\lambda}_2]&0              &  \ldots&0\\
0              & 0             &E[\widehat{\lambda}_3]&  \ldots &0\\
\ldots\\
0               &0              &0              &\ldots   &E[\widehat{\lambda}_p]  
\end{array}
\right)U^t\vec{a}+O\left(\frac{\sigma^2_j}{\sqrt{n}}\right).
\end{equation}
But now in the expression on the right side of the approximation above
only $U^t\vec{a}$ is random. What we can do now is to take the coefficients
of $U^t\vec{a}$ from another simulation and we would still
have the same approximation error size. We will simulate 
the eigenvectors of another matrix and write them in the basis of the 
eigenvectors of the estimated covariance matrix.    

More precisely, say that for another simulation $\vec{b}_j$ is the vector expressing
the $j$-th unitary eigenvector of the original covariance matrix 
$COV[\vec{Z}]$ expressed in the unitray basis of eigenvectors
of the estimated covariance matrix. Let
$b_{ij}$ be the $i$-th entry of the vector $\vec{b}_j$.
Taking the expectation of \eqref{lambda2} yields
\begin{equation}
\label{lambda3}
\lambda_j= E\left[\vec{b}_j^t 
\left(\begin{array}{cccccc}
E[\widehat{\lambda}_1]&0&0                            & \ldots&0\\
0              &E[\widehat{\lambda}_2]&0              &  \ldots&0\\
0              & 0             &\widehat{\lambda}_3&  \ldots &0\\
\ldots\\
0               &0             &0              &\ldots   &E[\widehat{\lambda}_p]  
\end{array}
\right)\vec{b}_j\right]+O\left(\frac{\sigma^2_j}{\sqrt{n}}\right)
=\sum_{i=1}^p E[b_{ij}^2]E[\widehat{\lambda}_i]+
O\left(\frac{\sigma^2_j}{\sqrt{n}}\right)
\end{equation}
where the last approximation above is up to a smaller order term.
Now, note that the coefficients $\widehat{\lambda}_i$ are known to us,
and since they are very close to their expectation
we assume that $E[\widehat{\lambda}_i]$ is known.

Now let $B$ be the matrix which has as $(i,j)$-entry
the expectation of the square coefficient $b^2_{ij}$.
This means that $B$ is the matrix who's $j$-column represents
the expected coefficient square of the representation of 
the $j$-th unitary eigenvector of the true covariance matrix
$COV[\vec{Z}]$ in the unitary basis of the eigenvectors
for the sample covariance matrix.   Hence 
\begin{equation}
\label{Bij}
B_{ij}=E[b_{ij}^2].
\end{equation}
Of course the matrix $B$ depends on the chosen spectrum 
$\vec{\lambda}$,
\[
\vec{\lambda}=(\lambda_1,\lambda_2,\ldots,\lambda_p)
\]
we took to simulate the matrix $B$.

Hence we write
\[
B=B(\vec{\lambda}).
\]
Recall our usual convention that $\lambda_1>\lambda_2>\ldots>\lambda_p$.

Thus, \eqref{lambda3} is an equation for every $j\in \{1,2,\ldots,p\}$.
We can write these $j$ equations in matrix/vector form as
\begin{equation}
\label{fundamental}
\vec{\lambda}=(E[\widehat{\lambda}_1],\ldots,E[\widehat{\lambda}_p])
\cdot B(\vec{\lambda})+\vec{\epsilon}(\vec{\lambda}),
\end{equation}
where $\vec{\epsilon}(\vec{\lambda})$ is a smaller order term,
which is to be considered a non-random function of the vector
$\vec{\lambda}$.  The right side of the above
equation \eqref{fundamental} can be viewed as a non-random function of the 
imput $\vec{\lambda}$. The key is that we can interpret \eqref{fundamental} 
as saying that the true spectrum $\vec{\lambda}$ is a fixed point
for that map given in \eqref{fundamental}. However,
the term $\vec{\epsilon}(\vec{\lambda}$ is a smaller
order term and not observable which we are going to neglect.

{\bf Our reconstruction algorithm 
for the true spectrum $\vec{\lambda}$ consists
of finding a fixed point for the map}
\begin{equation}
\label{themap}
\vec{\nu}\mapsto
(E[\widehat{\lambda}_1],\ldots,E[\widehat{\lambda}_p])\cdot B(\vec{\nu})
\end{equation}
where $\vec{\nu}=(\nu_1,\nu_2,\ldots,\nu_p)$ with 
$\nu_1>\nu_2>\ldots>\nu_p\geq 0$.  

Here, again $B(\vec{\nu})$ consists of the $p\times p$ matrix who's $ij$ entry
is the expectation of the square of the coefficient
of the $j$-th eigenvector of the original covariance
matrix written in the basis of the sample covariance matrix given that
the true spectrum is $\vec{\nu}$. In other words, given $\vec{\nu}$,
we simulate the data many times using $\vec{\nu}$ as the true
spectrum of the $COV[\vec{Z}]$. For every simulation we compute the matrix
representing the unitary eigenvectors of the original matrix expressed
in the basis of the unitary eigenvectors of the simulated sample
covariance matrix. Then, we take the average of these coefficients
square, over the different simulations. This gives us a good
approximation of $B(\vec{\nu})$.  

Now the key part is that we start with any spectrum $\lambda_0$
and iterate the map.  Otherwise stated we construct a sequence of spectrae 
\begin{equation}
\label{sequenceofnu}
\vec{\nu}_0,\vec{\nu}_1,\ldots.
\end{equation}
Each of the $\vec{\nu}_k$ is a positive $p$ dimensional vector such that

$$\vec{\nu}_1=(E[\widehat{\lambda}_1],\ldots,E[\widehat{\lambda}_p])B(\vec{\nu}_0)$$
and in general
$$\vec{\nu}_{k+1}=(E[\widehat{\lambda}_1],\ldots,E[\widehat{\lambda}_p])B(\vec{\nu}_k)$$

If the sequence \eqref{sequenceofnu} converges,
then it must converge to a fixed point. And since $\vec{\lambda}$ is a
fixed point, we assume that 
$$\vec{\lambda}\approx\lim_{k\rightarrow\infty}\vec{\nu}_k$$
up to smaller order error term.

\subsection{Free probability approach}
\label{freeprobapproach}
Assume $A$ and $B$ are two
$n\times n$ symmetric random matrices which we take to be independent of each other. 

Let $\mu$ and $\nu$ be two distributions on the real line. 
We choose the eigenvalues of $A$, such that the distribution of the eigenvalues converge to a 
random variable distributed according to $\mu$. Similarly assume that the eigenvalues of $B$ 
distribute in the limit according to $\nu$.  In addition, assume that the 
distribution of the matrices are invariant under conjugation by unitary matrices. 
We can construct such matrices as follows, we chose the direction of all the 
eigenvectors of $A$ and $B$ at random so that
these directions are equiprobably pointing into any direction
in space. For example take standard normal symmetric matrices and use the eigendecomposition
into the spectral part and the eigenvectors part.  The matrix $U$ formed 
by the eigenvectors is a unitaty matrix and we can form the matrix 
\[
 A=UDU'
\]
where the matrix $D$ is a diagonal matrix with the diagonal being the eigenvalues of $A$ 
which distribute in the limit according to $\mu$.  We do the same thing for $B$ and in addition we assume the choices 
independent of the choices of $A$.  
Thus we produced here $A$ and $B$ which are now independent of each other 
if we choose the eigenvalues and the eigenvectors independent of each other.

 Consider now the matrix product $AB$. Free probability theory
demonstrates that as the dimension of the matrices tends to infinity, 
the empirical distribution of the spectrum of $AB$ converges
in distribution \cite{voiculescu1991limit,voiculescu1992free} to the so called free multiplicative convolution 
of $\mu$ and $\nu$ and is denoted by $\mu\boxtimes \nu$. 
The best way to understand this is  via the so called $S$ transform introduced by 
Voiculescu and described in \cite{voiculescu1992free}.  This is very similar in spirit with the Fourier 
transform from the classical probability.  For any distribution $\mu$, it assigns a 
formal power series $S_{\mu}(z)$ which can be also understood as a complex analytic function 
in some region.  The main property of this transform is that 
\[
 S_{\mu\boxtimes \nu}(z)=S_{\mu}(z)S_{\nu}(z).
\]
 
Coming back to our problem, notice that 
$\vec{X}=(X_1,X_2,\ldots,X_p)$ is a normal vector with independent entries.  
Define now $\sigma^2_j:=VAR[X_j]=\lambda_j$.
We can write $\vec{X}$ as a diagonal matrix 
times a vector with standard normal entries as follows
\[
 \vec{X}=\vec{N}\Sigma
\]
where we define $\Sigma$ as the diagonal matrix with diagonal vector given 
by $\sigma_1,\sigma_2,\dots,\sigma_p$.

Now, all the rows of the matrix $X$ are i.i.d. and each is distributed
like the vector $\vec{X}$. Therefore we can write the matrix $X$ as
\[
X=N\cdot \Sigma
\]
where $N$ is a $n$ time $p$ matrix with all entries independent standard normal.
This in turn gives the estimated covariance matrix as
\[
\widehat{COV}[X]=\Sigma\cdot \frac{N^t\cdot N}{n}\cdot \Sigma.
\]

Since for any square matrices $C$ and $D$, the spectrum of $C\cdot D$ 
is equal to the spectrum of $D\cdot C$ we conclude that the spectrum of 
$\widehat{COV}[X]$ is equal to the spectrum of 
\begin{equation}
\label{rot}
\frac{1}{n}N^t\cdot N\cdot \Sigma^2
\end{equation}
If $p$ is large enough and $p/n$ stays approximately $r$, then the 
empirical distribution of \eqref{rot} is approximately the
free product of the spectrum of $\Sigma^2$ and the  matrix:
\begin{equation}
\label{marchenkopastur}
\frac{1}{n}N^t\cdot N.
\end{equation} 
The matrix \ref{marchenkopastur} is known
to converge to Marchenko-Pastur distribution given by 
\[
 \alpha(dx)=\frac{\sqrt{(\rho_{+}-x)(x-\rho_{-})}}{2\pi\rho x}1_{[\rho_{-},\rho_{+}]}dx
\]
where $\rho_{\pm}=(1\pm\sqrt{r})^2$
and $n$ and $p$ both go to infinity while $p/n$ is an approximately 
fixed proportion of $\rho$ which we take to be between $(0,1)$.

In the limit, the distribution of eigenvalues of $\widehat{COV}(X)$ becomes 
\[
\widehat{\mu}=\mu\boxtimes\alpha
\]
where $\mu$ is the limiting distribution of the spectrum of $\Sigma$.  
Thus, if we use the $S$ transform we get now 
\[
 S_{\widehat{\mu}}(z)=S_{\mu}(z)S_{\alpha}(z). 
\]
Recall that the $S$ transform is given by 
\begin{equation}\label{e:Stran}
S_\mu(z)= \frac{z+1}{z}m_{\mu}^{-1}(z)
\end{equation}
where the inverse is the formal power series inverse and 
we assume that $\int x\mu(dx)\ne0$ with
\[
m_{\mu}(z)=\sum_{m\ge1}\left(\int x^n\mu(dx)\right)z^n.
\]
Notice here that the inverse in \eqref{e:Stran} can also be made analytically by 
restricting the domain of $z$ to a proper domain in the complex plane.

From this we can invert the relation into 
\[
 S_{\mu}(z)=\frac{S_{\widehat{\mu}}(z)}{S_{\alpha}(z)}
\]
and thus 
\begin{equation}\label{e:fp:1}
 \mu=S^{-1}\left(\frac{S_{\widehat{\mu}}(z)}{S_{\alpha}(z)} \right)
\end{equation}
which is a way of recovering the distribution of the eigenvalues of $\Sigma$ in the limit.

Now, recall what our fundamental problem is.  
Given the spectrum of the estimated covariance matrix 
$\widehat{COV}[\vec{X}]$ reconstruct the spectrum of the original covariance
matrix $COV[\vec{X}]$, which is given by 
\[
\sigma^2_1>\sigma^2_2>\ldots>\sigma^2_p.
\]
If we consider that the estimated covariance matrix spectrum is
very close in distribution to the free product of the Marchenko-Pastur
and $\Sigma^2$ (because the dimension is big enough), then we can conclude that 
reconstructing the spectrum $\Sigma$ can be done using an 
approximation to \eqref{e:fp:1}.  Thus we have the following approximation  
\[
S_{\widehat{COV}[\vec{X}]}(z) =S_{\frac{N\cdot N}{n}\cdot\Sigma^2}(z)
\approx S_{\frac{N\cdot N}{n}}(z)\cdot S_{\Sigma^2}(z)
\]
which then leads to 
\[
\Sigma^2\approx 
S^{-1}\left( 
\frac{S_{\widehat{COV}[\vec{X}]}}{S_{\frac{N\cdot N}{n}}}
\right).
\]

This is a nice description in terms of free probability.  
On the other hand we still do not have a very good numerical 
technique to deal with this for finite reasonable large $n$ and $p$. 

This free probability description can also be phrased in terms of Bai-Silverstein approach
which essentially is the same characterization as \eqref{e:fp:1} only that the equation is 
writen explicitly.  This is for example treated in \cite{ledoit2015spectrum,ledoit2012nonlinear} 
which is inspired by \cite{el2008spectrum}.

\section{The analysis and the structure of the data}

\subsection{Simulations in the normal variables case}
Again take, 
\begin{equation}
\label{train}
\vec{\mu}_1,\vec{\mu}_2,\ldots,\vec{\mu}_p
\end{equation}
to be the principal components (eigenvectors) of the original
covariance matrix $COV[\vec{Z}]$ corresponding
to the eigenvalues
$$\lambda_1>\lambda_2>\ldots>\lambda_p.$$
We assume this eigenvectors to be unitary, that is they have
length $1$. The eigenvectors are apriory not known,
when we deal with real data. But they have nice
properties. We can express the vector $\vec{Z}$ 
in the basis of the principal components \eqref{train}.
Then the components are normal, independent
and and have variances corresponding to 
the eigenvalues of the covariance matrix.
More precisely, let $X_i$ be the $i$-th component
of $\vec{Z}$ expressed in that basis \eqref{train}. 
Hence
$$X_i=\vec{Z}\cdot \vec{\mu}_i.$$
Then, we have that $X_i$ is normal
and 
$$VAR[X_i]=\lambda_i$$
whilst $X_1,X_2,\ldots,X_p$ are independent
with expectation $0$.

Now, referring to the model in Example~\ref{e:ex2},
let $X_{ij}$ be the $j$-th component
of the portfolio on day $i$ expressed in the basis
of the eigenvectors \eqref{train}. Hence,
$$X_{ij}=\vec{\mu}_j\cdot(Z_1,Z_2,\ldots,Z_p),$$
Let
$$X=(X_{ij})_{ij}$$
be the $n$ times $p$ matrix
having as $i$-th, $j$-th entry $X_{ij}$
for all $i=1,2,\ldots,n$ and $j=1,2,\ldots,p$.

Hence, $X$ is the matrix whose $j$-the line gives
the return vector of the portfolio expressed in the basis
of the principal components. 

Now, note that $X^t\cdot X$ and $Z^t \cdot Z$
have the same spectrum since one can be obtained from
the other by an orthonormal change of basis since the 
principal components are orthonormal. 
In other words, if we are just interested in the spectrum
and how it changes when we estimate it, we can resimulate
$X$ instead of $Z$. Recall that $\vec{X}=(X_1,X_2,\ldots,X_p)$ has independent
normal components with expectation $0$. Note also that the rows
of $X$ are independent copies of $\vec{X}$. Thus,  to simulate the
estimated eigenvalues $\widehat{\lambda}_1,\widehat{\lambda}_2,\ldots,\widehat{\lambda}_p$
we are going to simulate $n$ independent copies of $\vec{X}$.

That is we simulate $\vec{X}_i=(X_{i1},X_{i2},\ldots,X_{ip})$
independently for $i=1,2,\ldots,n$. For each $i$, we let
$X_{i1},X_{i2},\ldots,X_{ip}$ be independent normal variables
with $VAR[X_{ij}]=\lambda_j$ and expectation. In this way we get the matrix
$$X=(X_{ij})_{ij}.$$
This matrix has the same distribution
as the stock matrix $Z$ expressed in the basis of the principal components.
Now, we take the estimated covariance matrix
$$\widehat{COV}[\vec{X}]:=\frac{X^t\cdot X}{n}$$
and this is the estimated covariance matrix
for the random vector $\vec{X}$. Note that the covariance matrix
of $\vec{X}$ is a diagonal matrix, since the entries
of $\vec{X}$ are independent. Furthermore, the estimated
covariance matrix 
$\widehat{COV}[\vec{X}]$ is simply the 
estimated covariance matrix $\widehat{COV}[\vec{Z}]$
expressed in the basis of the principal components
$\vec{\mu}_1,\ldots,\vec{\mu}_p$. So, the two 
estimated covariance matrices have the same spectrum.
Hence, to study, the spectrum of  $\widehat{COV}[\vec{Z}]$,
we can (and will) simply simulate $\widehat{COV}[\vec{X}]$,
since it has the same spectrum.

\subsection{Use of recovering the spectrum with normal data
and other data}
There are several useful things we can do in data-analysis
 once we are given the   spectrum
of $E[Z^tZ]$, (rather
than just the  spectrum of $Z^tZ$). For example, as one can see
in Figure \ref{1000vs2000}, the ``elbow'' is visible in
the reconstructed spectrum but not in the sample covariance matrix
spectrum. Another important issue is whether in the true spectrum
there are a lot of eigenvalues close to $0$ or not. In
stock data, there are no eigenvalues close to $0$ 
in the spectrum of $E[Z^tZ]$, but there will be in the spectrum
of $Z^tZ$ if the sample size is not quite a bit bigger than the vector
size. To understand why with stocks there is no
eigenvalues close to $0$, simply consider in our previous example 
in Section~\ref{principalcomponents} and notice that the is always the term
$\epsilon_i$. There is no stock which can be almost
entirely predicted by the index of the sector there is always
a sizable firm-specific part $\epsilon_i$. And indeed, when
we reconstruct the spectrum of the stocks correlation matrix, 
the smallest eigenvalues are around $0.3$ and not close to $0$.
Now, when we consider financial futures the situation is very different. 
There are several futures which are highly correlated leading to eigenvalues
in the correlation matrix close to $0$. For example, silver and gold futures
are highly correlated. Since our reconstruction method is more precise than
others, allow much better distinction between small eigenvalues
in $Z^tZ$ due to error because of small sample,  and ``real'' close to $0$
eigenvalues in $E[Z^tZ]$.

Let us give  two more applications, of what to do with the 
spectrum of the true covariance matrix or kernel.  
\begin{enumerate}
\item{}If the  data is multivariate normal, 
we can simulate the whole situation using the true spectrum.
Indeed, for simulating the relative position of the estimated 
principal components with respect to the true principal components,
there is no other parameters needed than the spectrum of
the true covariance matrix, (provided we have normal data).
We can see how the estimated principal components behave with 
respect to the true principal components, by simulation
using the reconstructed ``true'' spectrum. That is if a practitioner
asks ''what do these principal components represent'' one can do a simulation
with  the ``true'' reconstructed spectrum and find out the relative position
of the principal components of the $Z^tZ$ versus the principal components
of $E[Z^tZ]$. We can find out exactly how these estimated principal components
lie in the space with respect to the true principal component. If they are
too far away from all eigenvectors with close eigenvalues, one can 
immediately see that the estimated principal components
may have  little meaning. Similarly, in Linear Discriminant Analysis (LDA), 
assuming the data to be normal, we can resimulate the data using 
the true spectrum and find out which way of correcting 
the estimated covariance matrix works best when $n$ is not much larger than
$p$. Indeed,  in LDA, we have that the classification is made based on
a formula involving the inverse of the estimated covariance matrix. 
But when $n$ is not much
bigger than $p$, then there will be a large amount of eigenvalues
of the estimated covariance matrix close to $0$, which in reality 
are likely to be very different from $0$. This will result
in an enormous error as far as the inverse of the 
covariance matrix is concerned. 

\item{} Detailed analysis of the underlying structure
of real life data is possible with  
the spectrum of $E[Z^tZ]$, rather than the spectrum of $Z^tZ$.
For example we can take subsamples from the data and see how the spectrum
behaves compared to the spectrum of the whole data. But, for this we
need the spectrum of the underlying ``true'' covariance matrix, rather
than the spectrum of the sample covariance matrix. We give an example
of such an analysis below. 
\end{enumerate}
 
{\footnotesize {\bf Example of analysis of the spectrum using subsamples.}
We take the stocks daily-return-data of $800$ stocks and 
$2000$ days to find numerically 
the spectrum of the true underlying covariance matrix
using our algorithm. Then, we do the same 
for our data-set restricted to $400$ stocks.
That means that we pick $400$ stocks at random among our $800$
and then restrict the data-set to the selected stocks.
Finally, we restrict the data set to $200$ randomly chosen stocks.
For all three numbers of stocks,
 we calculate the ``true underlying spectrae'' and plot
them joint onto the same plot in order to compare the three different
empirical distribution of the spectrum.
We adjust the indexes, so that each plot has the same
$x$-coordinate-width. (The rescaling is done so that the $800$-hundredths eigenvalue of the $800$-stock data,
is represented at the same $x$-coordinate then the 
$400$-th eigenvalue of the data restricted to $400$ stocks.
So, we do this by plotting the $800$ eigenvalues against
$1,2,\ldots,800$ and plotting the $400$ eigenvalues
against $2,4,6,\ldots,800$.
 
 Let $\lambda_i^{800}$ designate the $i$-th
eigenvalue of the covariance matrix with the full data
of $800$ stocks. Let $\lambda_i^{400}$ be the $i$-th eigenvalue
of the data restricted to $400$ stocks chose at random.
We plot in black the points
$$ (1,\lambda^{800}_1),(2,\lambda^{800}_2),
(3,\lambda^{800}_3),\ldots,(800,\lambda^{800}_{800}),$$
then in red color we plot the following points
$$(2,\lambda^{400}_1),(4,\lambda^{400}_2),(6,\lambda^{400}_3),
\ldots,(800,\lambda^{400}_{400}).$$
and finally, we plot
$$(4,\lambda^{200}_1),(4,\lambda^{200}_2),(14,\lambda^{200}_3),
\ldots,(800,\lambda^{200}_{200}),$$
where $\lambda^{200}_i$ represents the reconstructed
eigenvalue number $i$ of the data set with $200$ stocks.

We see a surprising thing these curves look similar and have 
there ``elbow'' approximately
at the same place, when we rescale them properly, except maybe for
the 15 or so biggest eigenvalues. For the biggest eigenvalues we have
linear growth, which is very different from giving the 
same ``empirical distribution function''.

The result is seen in Figure \ref{compare_different_p} below:
\begin{figure}[!ht]
\caption{Comparing reconstructed spectrum's empirical distribution
for different numbers p of stocks}
\label{compare_different_p}
\includegraphics[width=90mm]{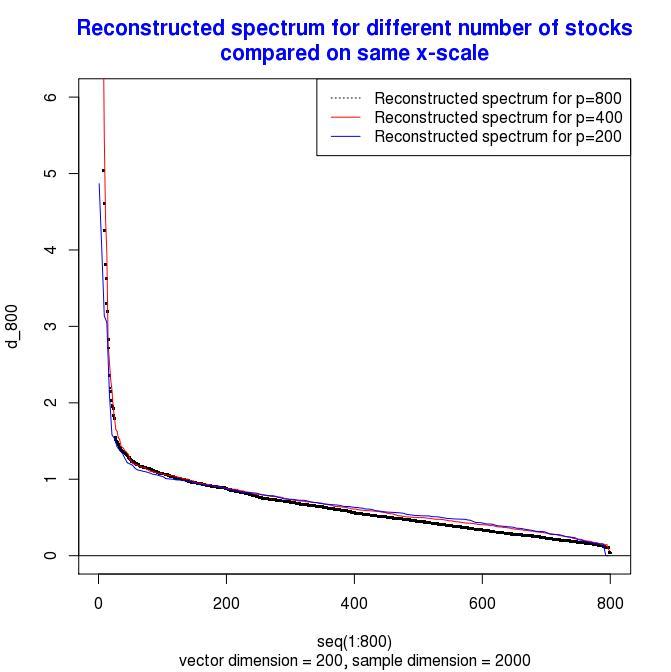}
\end{figure}

 Now Figure \ref{compare_different_p} contradicts the simple model  
presented in section \ref{principalcomponents}.
The reason is this.   With that model for stocks
as being a linear combination of a sector indexes plus an
independent firm specific term $\epsilon_i$, we would have 
a fixed number of factors (mainly sectors) which lead to 
eigenvalues of order $O(p)$. The number of these factors does 
not change whether we have $400$ or $800$ stocks, since the 
stocks are chosen at random.  In fact, in our actual data there are $11$ sectors.
These sectors and a few other factors are responsible
for the biggest eigenvalues.
(For details about why these eigenvalues are of order $O(p)$
consider our previous example of covariance matrix 
of stocks $Z,Y,Z,W$ given at the very  end of Section \ref{sectionintro}). But, then
if the elbow which is seen in the spectrum 
in Figure \ref{compare_different_p} would correspond
to the transition between the order $O(p)$ eigenvalues and the eigenvalues
due to the diagonal, then it should happen at the same index for the $800$
sector data and for the $400$ sector data. By this we mean the same
index in absolute terms and not in relative terms. Instead,
the elbow is observed at the same relative position, i.e.
when we represent the two spectrae above the same $x$-coordinate
interval, so as to compare the empirical distribution of the eigenvalues
of the restricted and the non-restricted data-sets.
The second thing is that the values where the elbow happen are 
eigenvalues above one. This also contradicts a simple model 
containing only  a number $O(1)$ of order $O(p)$ factors
 plus the  firm independent terms and nothing else. Indeed,
according to the model for stocks presented
in Section \ref{principalcomponents}, apart from the big factors which lead
to $O(1)$ eigenvalues of order $O(p)$, the rest of the spectrum
should be close to the set $\{VAR[\epsilon_i]|i\in 1,\ldots,p\}$.
However, since we rescale the stocks so that they have variance $1$,
that is we present spectrum of correlation matrix, we have
that the variances $VAR[\epsilon]_i<1$ if the $\epsilon_i$'s
are uncorrelated to the factors. This would imply that
 the eigenvalues which are not order $O(p)$,
should be strictly less than $1$.

Let us next see a table with  the ratio between the first
30 eigenvalues of the data with $800$ stocks divided by 
the eigenvalues of the data with $800$ stocks. We present
this ratio for the $30$ biggest eigenvalues.
The eigenvalues which depend on a factor present in a large part of stocks
should grow linearly in $p$. (See subsection \ref{our_modelling} for details
as to why).  In the table below we see that indeed the biggest eigenvalues
are proportional to $p$. So, when we double $p$, these eigenvalues
should also approximately double. Hence, the biggest eigenvalue should
roughly double as we go from the data with $400$stocks to the data with
$800$-stocks. This is indeed the case for the first $14$ eigenvalues
as can be seen in the next table: 

$$
\begin{array}{c|c|c|c|c|c|c|c|c|c| c|c|c|c|c|c}
\lambda^{800}_i/\lambda^{400}_i& 2.1& 2.0& 2.0& 2.1& 1.8& 1.9& 2.0& 1.8& 1.7& 1.9& 1.8& 1.9& 1.8& 1.9& 1.7\\\hline
i&1&2&3&4&5&6&7&8&9&10&11&12&13&14&15        
\end{array}
$$

and

$$
\begin{array}{c|c|c|c|c|c|c|c|c|c| c|c|c|c|c|c}
\lambda^{800}_i/\lambda^{400}_i&1.7& 1.6& 1.6& 1.5&
 1.4& 1.4& 1.4& 1.4& 1.4& 1.4& 1.2& 1.2& 1.3& 1.3& 1.3\\\hline
i&16&17&18&19&20&21&22&23&24&25&26&27&28&29&30        
\end{array}
$$

After eigenvalue $14$, it seems that the linear growth quickly stops
as the ratio becomes much smaller than $2$. Now, the elbow
in Figure \ref{compare_different_p} for the data with $800$ stocks
happen much after the $14$-th eigenvalue. (At least the lower part of the 
elbow). This implies that in our data there is more than just a number $O(1)$
of factors present in a large part of the stocks and the
firm specific terms, additional to that there must be an order
$O(p)$ of ``middle size'' eigenvalues coming from small number
of stocks interactions.

We presented this partial analysis here for one reason. This would
not be possible without our reconstruction formula for the spectrum!
Indeed, if we do not know the spectrum of the ``true'' underlying
covariance matrix, we can not observe  useful facts like for example
which eigenvalues grow linearly. Indeed, if you take for your analysis
sample covariance matrices spectra, you would not necessarily see
that you are dealing with the same spectral distribution for
most of the eigenvalues except those of order $O(p)$.}

\subsection{Big data}
Assume that you have a correlation matrix of thousand by thousand.
Typically finding all the eigenvalues and eigenvectors goes quite fast,
on a personal laptop using $R$. Now, for $10000$ times $10000$
general correlation matrix (not sparse)  might be more difficult
and we may need to run it for a longer period. Often time, 
in big data you may have a much bigger  dimension than that 
of $10000$ times $10000$,  for example a million times a million. 
Then, with a classical algorithm  and a regular laptop, we
can not find all eigenvectors and eigenvalues in short time.
Thus,  there are some approximation algorithms which can handle 
much bigger data for finding eigenvectors and eigenvalues.
But from experience we know that in certain
situations they may not work. However, imagine that by one way
or the other for a correlation matrix of a million by a million,
we manage to get the sample spectrum. Then, our methods can 
easily ``reconstruct the true underlying spectrum''.
How? Well first note that our method with the asymptotic formula
for the high dimensional case (formula \ref{approxlambda}), requires only $O(p^2)$
calculation steps and what is more important only order $O(p)$ working
memory. The ``curse'' with finding eigenvectors for high dimensional
correlation matrix when $n$ is of order $O(p)$ is that
typically the number of calculation steps is  at least of 
order $O(p^3)$ and may be larger depending on 
the size of the spectral gaps. This assuming that the matrix does not
have a special structure which makes it easier to compute.
So, once we have the spectrum of the sample covariance matrix
given to us, our method based on the asymptotic multidimensional formula
is very easy to compute even with giant data-sets. What about
our fixed point method? Note that method also only requires
as input the spectrum of the sample covariance matrix,
as well as the two numbers $n$ and $p$. However, that method
tries to find a fixed point for a map associated with the second moments
of the coefficients of the matrix of eigenvectors of the true
covariance matrix expressed in the basis of principal components
of the sample covariance, which has to be simulated. More precisely,
 we have to simulate the data 
using a first guess of what the spectrum might be in reality. 
Then for that simulated data, we have to compute all the 
eigenvectors. And since we want the eigenvectors coefficients 
second moments, we have to repeat this a few times.
This is going to be
costly in computational time.  This might not be achievable
in very high dimension. But, the good news is that we can scale things
down. That is if you are given $n$ and $p$ simply solve the problem
for much smaller data set where the ratio $n/p$ remains the same
and where you take approximately the same empirical distribution
of the sample covariance spectrum. Once, you find in this way
the reconstructed spectrum you scale it up and you get almost
the exact answer. In the example below this is what we did:
we divided $n$ and $p$ by a factor $5$. The big data was $2000$ days 
and $800$ stocks. Then, we simulate a data-set which is $5$ times
in the number of samples and the size of vector. That is we simulate
a data set which would correspond to $160$ stocks and $400$ days.
We took for this a five times dimension reduced spectrum from the full
sample covariance. That is we took the full data sample spectrum.
Then we produced a spectrum on a dimension five times less, so that
the new spectrum would have approximately the same empirical
distribution. (Notice that we did not scale down the real life data, but
only the sample spectrum). Then, we do simulations for our reconstruction
on the reduced dimensional space. The resulting reconstruction is then presented
on the same scale as the higher dimension reconstruction.
The result can be seen in the next figure below
(see Figure \ref{compare_different_p}, where we compare the reconstructed on the 
full data, versus the reconstructed on the scaled down data by a factor $5$.
 We see that the scaled down reconstruction does an almost perfect job
in reconstructing the ``true higher dimensional spectrum''.
 Again, here we did not scale down the data, but
only the reconstruction process! This means that from
the given full spectrum of the sample covariance matrix,
we scale down by taking a spectrum with basically approximately same 
empirical distribution but vectors on a five times smaller dimension.
Then, we did the fixed point eigenvector reconstruction on
that smaller space. The result is shown again in the higher
dimensional plot, but showing the obtained empirical spectrum
on the right scale so that it can be compared to the reconstructed
spectrum at full scale. And indeed we got perfect match. Note
that we had to adjust the 5 biggest eigenvalues a little bit,
to obtain the result.
\begin{figure}[!ht]
\caption{Reconstruction from full sample spectra, scaled down to a lower
dimensional space, a simulations for reconstruction done in lower dimensional
space}
\label{reconstruction_with_scaling}
\includegraphics[width=90mm]{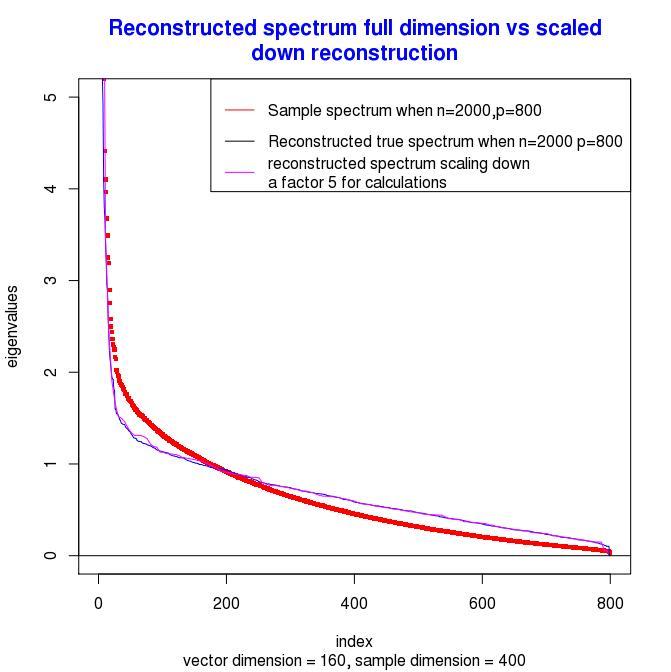}
\end{figure}

\section{Testing the different methods numerically}
In real data we have always observed spectrae which look continuous
and highly convex. But in the  model in  section \ref{principalcomponents},
the stocks were written as linear combination 
of factors and the part which comes only from
the firm. Thus in theory according to  such a model, we would 
have  a small number of order $O(1)$ of eigenvalues of order $O(p)$. The 
rest of the eigenvalues would be of order $O(1)$.
In reality, this is never what we have seen in real-life data.
Instead, we have a continuum going from the biggest eigenvalues
 to the eigenvalues of order $O(1)$. Indeed there is often,
one eigenvalue (or a very small number) which is separated from
the rest. After this largest eigenvalue (or couple of them)
we don't go directly to something which is order $O(1)$ and does
not grow, but instead have a continuum of eigenvalues among which some 
grow as we increase $p$ and hence are not of order $O(1)$. 
Hence, to evaluate our spectrum reconstructing methods, we are primarily
interested in such spectrum with the properties we encounter mostly in real data, 
that is, ``continous functions except mabye for the biggest eigenvalue
(or maybe an very small number of eigenvalues)'' and the rest of the eigenvalues
having the function $s\mapsto \lambda_s$ looking continuous convex
and with derivative going to $-\infty$ as the index goes to $0$.  
Recall that the eigenvalues of $E[Z^tZ]$ are denoted 
by $\lambda_1>\lambda_2>\ldots>\lambda_o$.

In our subsequent simulations, we often take a relatively small 
dimension with $p=100$ for example. Why?  To make simulations
quicker and because we have seen that even if you blow up things from there
things do not change too much. In Figure \ref{blowup1000} we see what happens
when we take $n$ and $p$ both $100$ times bigger, but we take about
the same empirical distribution for the spectrum of the original covariance
matrix $COV[\vec{X}]$.  The empirical distribution of the estimated covariance
matrix $\widehat{COV}[\vec{X}]$ barely changes. Thus, to test the different recovering
methods, a low dimension like $p=100$ should be enough in practice
when dealing with synthetic data. With real life data, things are very different,
and we also need to test our method with large $p$. Indeed,
when $p$ becomes larger, the real spectrum of the underlying
covariance matrix does not have approximately apriori the same
distribution, since the biggest eigenvalues are supposed to
grow linearly in $p$!

Another important thing, would be to distinguish between those
spaectrae (of underlying covariance matrix)
 which are strongly bounded away from $0$ (like $\lambda_s\geq 0.3$ for
example in stocks) and the others like in futures, where
there a close to $0$ eigenvalues, reflecting the possibility
of arbitrage. Indeed if the true spectrum is bounded away 
from zero say for example by $0.3$, then
a good measure for comparing how well the different 
reconstruction methods work, would be to take the maximum
of relative error:

$$\max_{s=1,\ldots,p}\frac{|\lambda_s-\widehat{\lambda}_s^M|}{\lambda_s},$$
where $\widehat{\lambda}_s^M$ is the estimate of the $s$-th eigenvalue
$\lambda_s$
using the estimation method $M$.  If there are eigenvalues $\lambda_s$
close to zero we need to use another method to access the precision
of the reconstruction method as in that case the relative error for the values
$\lambda_s$ close to $0$ is meaningless. So, for that case we will
take the maximum of the relative error only for those $\lambda_s$
which are bigger than $0.2$. For those which are smaller we take the
maximum absolute error.

\subsection{Testing our methods on synthetic multivariate normal data}
 \label{testingdata}

We simulate multivariate normal data using various
spectrum: linear, two different linear segments combined, step function and spectrum from stock data. In particular, we compare our eigenvector method with Quest provided by the free probability approach \cite{ledoit2012}.

\begin{figure}[!ht]
			\begin{subfigure}[t]{.5\textwidth}
			\centering
			\includegraphics[width=.9\linewidth]{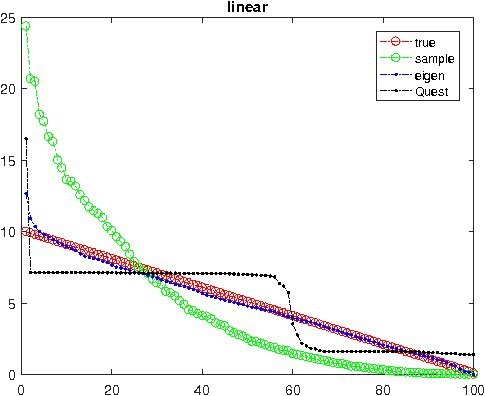}
		\end{subfigure}
		\begin{subfigure}[t]{.5\textwidth}
			\centering
			\includegraphics[width=.9\linewidth]{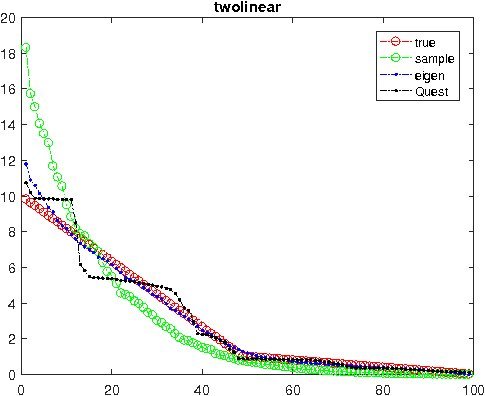}
		\end{subfigure}	
		\begin{subfigure}[t]{.5\textwidth}
		\centering
		\includegraphics[width=.9\linewidth]{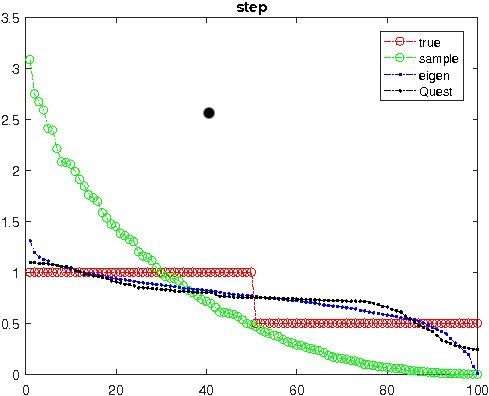}
	    \end{subfigure}
    	\begin{subfigure}[t]{.5\textwidth}
    	\centering
    	\includegraphics[width=.9\linewidth]{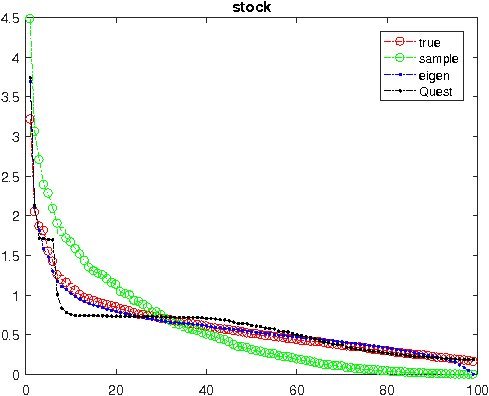}
        \end{subfigure}
	\caption{}
	\label{fig:methods_compare}
   \end{figure}

In the majority of those cases Figure \ref{fig:methods_compare}, eigenvector approach outperform others.  There is another significant advantage of eigenvector approach.  That is its insensitivity to the sample spectrum, where Quest varies quite significantly even when sample spectrum exhibits a small change. In other words, when a different sample matrix is given, the recovery of our eigenvector method does not fluctuate, which means the variance is small.

For a totally flat spectrum the first method does not work well: This is not too surprising
since our first method is based on comparing the $i$-th eigenvector
of the sample covariance matrix with the $i$-th eigenvector
of the true covariance. In the flat spectrum case, all eigenvectors
of the true covariance matrix are identical.
This can be seen in figure \ref{constantspectrum} below.
\begin{figure}[!ht]
\caption{First method fails, second not very precise for constant spectrum.}
\label{constantspectrum}
\includegraphics[width=90mm]{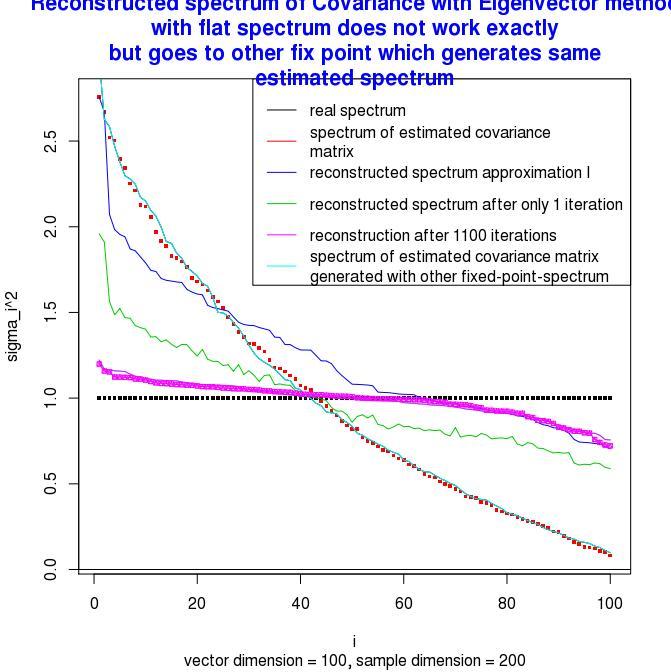}
\end{figure}  

\begin{figure}[!ht]
\caption{By scaling down and then reconstructing before blowing up again, we manage to retrieve the spectrum for big data.}
\label{Reconstructed_by_scaling_down}
\includegraphics[width=90mm]{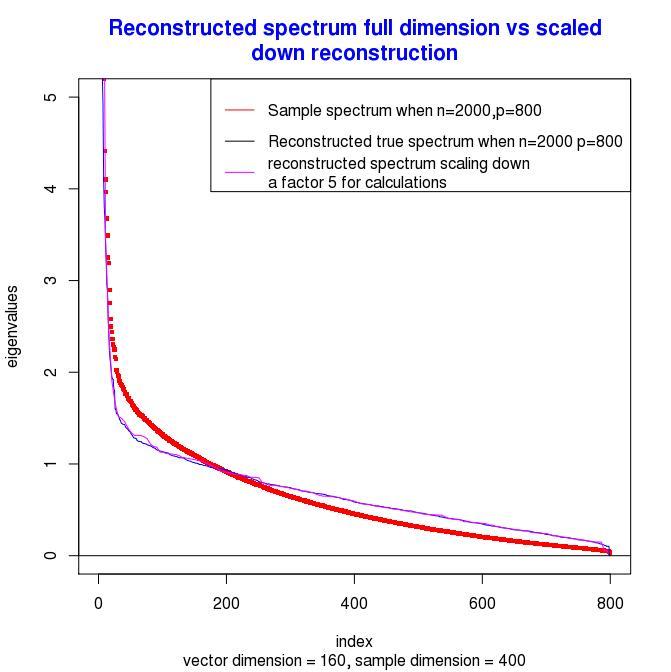}
\end{figure}

Also, often by scaling down and then retrieving the spectrum for the smaller
data, set we can retrieve the true spectrum. This is shown in Figure 
\ref{Reconstructed_by_scaling_down}.

\end{document}